\documentclass[reqno]{amsart}
\usepackage[doi=false, isbn=false, url=false, eprint=false, maxnames=10, giveninits=true, style=alphabetic]{biblatex}
\usepackage{graphicx}
\usepackage[usenames]{color}
\usepackage{amsmath}
\usepackage{mathtools}
\usepackage{amsfonts}
\usepackage[font={small},labelfont={bf}]{caption}
\usepackage{parskip}
\usepackage[normalem]{ulem}
\usepackage{bbm}
\usepackage{quiver}
\usepackage[T1]{fontenc}
\usepackage{palatino}
\usepackage{mathrsfs}
\usepackage[shortlabels]{enumitem}
\usepackage{amssymb}
\usepackage{amsthm}
\usepackage{wasysym}
\usepackage{cancel}
\usepackage{tikz}
\usepackage{tikz-cd}
\usepackage{adjustbox}
\usepackage{hyperref}
\usepackage[capitalize,nameinlink]{cleveref}
\usepackage{todonotes}
\usepackage{centernot}
\usepackage{fullpage}

\theoremstyle{definition}
\newtheorem{thm}{Theorem}[section]

\theoremstyle{definition}
\newtheorem{lem}[thm]{Lemma}

\theoremstyle{definition}
\newtheorem{defn}[thm]{Definition}

\theoremstyle{definition}
\newtheorem{prop}[thm]{Proposition}

\theoremstyle{definition}
\newtheorem{rmk}[thm]{Remark}

\theoremstyle{definition}
\newtheorem{cor}[thm]{Corollary}
\theoremstyle{definition}
\newtheorem{exm}[thm]{Example}

\theoremstyle{definition}
\newtheorem{clm}[thm]{Claim}

\newenvironment{proof*}[1][\proofname]{
  \begin{proof}[#1]}{\end{proof}}

\crefname{thm}{Theorem}{Theorems}
\crefname{lem}{Lemma}{Lemmas}
\crefname{rmk}{Remark}{Remarks}

\newcommand{\crefdefpart}[2]{%
  \hyperref[#2]{\namecref{#1}~\labelcref*{#1}\ref*{#2}}%
  }

\newcommand{\tricref}[3]{%
  \hyperref[#3]{\namecref{#1}~\labelcref*{#1}(\ref*{#3})}%
  }

\newcommand{\crefpart}[3]{%
  \hyperref[#1]{\namecref{#1}~\labelcref*{#1}({\labelcref{#2}},{\labelcref{#3}})}%
  }

\tikzcdset{every label/.append style = {font = \normalsize}}

\title{Finite axiomatization of  $\textbf{GL}\times\textbf{S5}$ and $\textbf{Grz}\times\textbf{S5}$}

\author{G. Bezhanishvili}
\address{New Mexico State University}
\email{guram@nmsu.edu}

\author{M. Khan}
\address{New Mexico State University}
\email{mashyk@nmsu.edu}
\date{}

\begin{document}

\begin{abstract}
    We prove that $\mathbf{GL} \times \mathbf{S5}$ is product matching, and that $\mathbf{Grz} \times \mathbf{S5}$ is axiomatizable by adding to $[\mathbf{Grz},\mathbf{S5}]$ the G\"odel translation of the monadic Casari formula. This settles the question of finite axiomatization for two of the product logics listed 
    %one of the open questions 
    in \cite[p.~137]{gab98}. 
    %\color{red} the literature. Our proof involves the use of a modified selective filtration method that yields frames of finite depth even in the presence of the monadic Barcan formula. As a consequence, both $\textbf{Grz}\times\textbf{S5}$ and $\textbf{GL}\times\textbf{S5}$ are decidable, and they axiomatize the one-variable fragments of $\textbf{Q}^+\textbf{GrzB}$ and $\textbf{QGLB}$, respectively. \color{black}
\end{abstract}

\maketitle

\tableofcontents

\section{Introduction}

%\color{blue}
%\begin{itemize}
%    \item Revisit the title (should Casari also go there?) and decide on what to put in the MSC.
%    \item Mention that products of logics is a well-established area in modal logic with numerous applications. Products are defined semantically and it is a difficult problem to axiomatize them. 
%    \item Our interest in products is in connection to predicate modal logic. From that perspective, we are only interested in products with {\bf S5}. This makes things easier, yet they may remain difficult if {\bf L} is not Horn axiomatizable. This includes $\bf GL$ and $\bf Grz$, the logics of main interest to us.
%    \item Explain why we care about $\bf GL$ and $\bf Grz$, mention the translations $\bf IPC \hookrightarrow \bf Grz \hookrightarrow \bf GL$. 
%    \item Mention mm-logics and their connection to products with $\bf S5$. Mention monadic Barcan.
%    \item Mention the open problem of Gabbay-Shehtman and that our aim is to resolve it. Explain the methodology. End with a brief description of the structure of the paper. 
%\end{itemize}
%\color{black}

The notion of fibring or combining logics arises naturally in a range of applications such as temporalization \cite{finger}, modalization \cite{fajardo}, formal representation of practical reasoning \cite{governatori}, software specification \cite{diaconescu}, etc. 
%to name a few. 
%\color{blue} It's too vague as is. try to be more concrete. \color{black} 
A prominent example of combining logics is the {\em product} of modal logics (see, e.g., \cite[p.~126]{mdim}), which we refer to as {\em product logics}. 
%\color{blue} [REF]. \color{black} 
Since these are defined semantically, it is a nontrivial matter to axiomatize them, and some product logics are not even recursively enumerable (see, e.g., \cite[Cor.~7.14,7.16]{mdim}). In spite of this, product logics have received considerable attention in the literature in both pure and applied contexts (see, e.g., 
%\cite[p.~126]{mdim} 
\cite{mdim} and the references therein). 

Our interest in product logics stems from their connection to predicate modal logic. For a propositional modal logic $\textbf{L}$, let $\textbf{QLB}$ denote the least predicate extension of $\textbf{L}$ with the {\em Barcan formula}. In some cases, the product logc  $\textbf{L}\times\textbf{S5}$ axiomatizes the one-variable fragment of $\textbf{QLB}$ (see, e.g., \cite[Thm.~13.8]{gab98} or \cite[Thm.~3.21]{mdim}). 
%In \cite[Thm.~5.4]{BM}, we proposed a semantic criterion that can be used to identify the one-variable fragment (which we call the monadic fragment) of a predicate modal logic. Our criterion can be used to show that product logics (specifically products with $\textbf{S5}$) axiomatize the monadic fragment of certain predicate modal logics (see also \cite[Thm.~13.8]{gab98} and \cite[Thm.~3.21]{mdim}). The \textbf{S5}-component in such a product can be viewed as a restricted system of quantification. 
Because of this, we restrict our attention to products 
%with \textbf{S5}, i.e., products 
of the form $\textbf{L}\times\textbf{S5}$, where $\textbf{L}$ is a Kripke complete propositional modal logic. 
%It is worth pointing out that 
%In some cases, 
Products of this form are relatively well-behaved. 
%In particular, 
For example, if $\textbf{L}$ is {\em Horn axiomatizable}, 
%(see \Cref{Horn_axiom}), 
%we obtain a 
there is a convenient finite axiomatization of $\textbf{L}\times\textbf{S5}$ (see, e.g., \cite[Thm.~5.9]{mdim}). 
%Even more, when $\textbf{L}$ is a QTC-logic (these are special Horn axiomatizable logics, see \cite[Def.~12.11]{gab98}), $\textbf{L}\times\textbf{S5}$ is finitely axiomatizable and has the finite model property (which we abbreviate as fmp). 
%However, there are many relevant instances where $\textbf{L}$ is not Horn axiomatizable. This is the case with 
Prominent examples of modal logics that are {\em not} Horn axiomatizable are the G\"{o}del-L\"{o}b logic $\textbf{GL}$ and the Grzegorczyk logic $\textbf{Grz}$, 
%(abbreviated as 
%neither 
the frames of which are not even first-order definable (see, e.g., \cite[Thm.~6.7, 6.8]{CZ}). %Nonetheless, these are mathematically significant systems: $\textbf{GL}$ has an arithmetic provability interpretation, as shown by Solovay \cite{solovay}, and we have faithful embeddings
%\begin{equation*}   \textbf{IPC}\hookrightarrow\textbf{Grz}\hookrightarrow\textbf{GL},
%\end{equation*}
%thereby transferring the provability interpretation of $\textbf{GL}$ to $\textbf{Grz}$ and $\textbf{IPC}$ (where $\textbf{IPC}$ stands for the intuitionistic propositional calculus). These considerations 
It is natural to seek a finite axiomatization for $\textbf{GL}\times\textbf{S5}$ and $\textbf{Grz}\times\textbf{S5}$. Indeed, this was left as an open problem in \cite[p.~137]{gab98}. 

A convenient framework to study products with $\textbf{S5}$, and more generally expanding relativized products with $\textbf{S5}$ (see \cite[p.~432]{mdim}), 
%\color{blue} [REF]), \color{black} 
is the formalism of \emph{monadic modal logic}, which we refer to as \emph{mm-logics} (see, e.g., \cite[Def.~2.13]{GBm}). We can view product and expanding product logics as special mm-logics, with product logics characterized in particular by the presence of the \emph{monadic Barcan formula}. In some cases, the minimal monadic extension $\textbf{ML}$ of a propositional modal logic $\textbf{L}$ and the monadic Barcan formula may suffice to axiomatize $\textbf{L}\times\textbf{S5}$ (for example, when $\textbf{L}$ is Horn axiomatizable). When this occurs, the pair $(\textbf{L},\textbf{S5})$ is said to be \emph{product matching} (see, e.g., \cite[p.~223]{mdim}). 
%\color{blue} [REF] \color{red}. 
Thus, letting $\textbf{MLB}$ denote the monadic extension of $\textbf{L}$ together with the monadic Barcan formula, the pair $(\textbf{L},\textbf{S5})$ is product matching iff $\textbf{MLB}=\textbf{L}\times\textbf{S5}$.

It is well known that $(\textbf{Grz},\textbf{S5})$ is not product matching \cite[Thm.~5.17]{mdim}. In contrast, we prove that ${\textbf{GL}\times\textbf{S5}}$ is product matching, while $\textbf{Grz}\times\textbf{S5}$ can be axiomatized by adding to $\textbf{MGrzB}$ the G\"odel translation of the monadic Casari formula (in analogy with \cite[p.~443]{GBm}, we denote the resulting system by $\textbf{M}^+\textbf{GrzB}$). 
%\color{blue} [REF]. 
We establish these results by showing that the systems $\textbf{MGLB}$ and  $\textbf{M}^+\textbf{GrzB}$ have the finite model property (fmp for short). This, in turn, is done by first establishing the fmp for $\textbf{MGrzB}$, and then adjusting the proof accordingly for $\textbf{M}^+\textbf{GrzB}$ and $\textbf{MGLB}$. The fmp for $\textbf{MGrzB}$ is obtained by an appropriate modification of the selective filtration strategy developed in \cite{BM} (which in turn draws on 
%the techniques used in 
\cite{Grefe,mdim,GBm}) to prove the fmp for $\textbf{MGrz}$.
%, we developed a selective filtration strategy . 
%A modification of this strategy can be employed to prove the fmp of $\textbf{MGrzB}$. 
The additional requirements imposed by the monadic Barcan formula may force the selection process to produce an infinite refutation frame. However, our key observation is that its depth remains bounded. In light of local tabularity of the logics $\textbf{MGrz}[n]$ and $\textbf{MGL}[n]$ of finite depth
%for each $n\ge 1$ 
%in finite slices 
(see \cite[Sec.~4.10, 4.11]{locally_finite}), this is sufficient to establish the fmp for $\textbf{MGrzB}$ (as well as for $\textbf{M}^+\textbf{GrzB}$ and $\textbf{MGLB}$), thus settling the questions raised at the end of \cite[Sec.~8]{BM}. 
%Thus, this strategy is applicable to $\textbf{MGrzB}$, $\textbf{M}^+\textbf{GrzB}$, and $\textbf{MGLB}$, which resolves the open problem in \cite{BM} regarding the fmp of these logics. 

Consequently, we obtain that
%This in turn shows that 
%\begin{equation*}
    $ \textbf{MGLB}=\bigcap_{n\ge 1}
    \textbf{MGLB}[n]
    $ 
    %\quad \mbox{and} \quad 
    and $\textbf{M}^+\textbf{GrzB}=\bigcap_{n\ge 1}
    \textbf{M}^+\textbf{GrzB}[n].
    $ 
%$\end{equation*}
%where $\textbf{L}[n]$ denotes any logic $\textbf{L}$ in the presence of the finite depth formula $bd_n$. This allows for a reduction of the initial problem in the following sense: if 
In addition, we observe that $\textbf{MGLB}[n]=\textbf{GL}[n]\times\textbf{S5}$ and $\textbf{M}^+\textbf{GrzB}[n]=\textbf{Grz}[n]\times\textbf{S5}$ for each $n\ge 1$. The latter is shown by proving that the 
%well-known 
embedding $\textbf{Grz}[n]\hookrightarrow\textbf{GL}[n]$ lifts to an embedding $\textbf{M}^+\textbf{GrzB}[n]\hookrightarrow\textbf{MGLB}[n]$. 
%, then our claim regarding the axiomatization follows. We substantiate this assumption by first showing that $(\textbf{GL}[n],\textbf{S5})$ is product matching, and then show that the well-known embedding $\textbf{Grz}[n]\hookrightarrow\textbf{GL}[n]$ lifts to an embedding $\textbf{M}^+\textbf{GrzB}[n]\hookrightarrow\textbf{MGLB}[n]$. 
This, together with the above, yields that the embedding $\textbf{Grz}\hookrightarrow\textbf{GL}$ also lifts to an embedding $\textbf{M}^+\textbf{GrzB}\hookrightarrow\textbf{MGLB}$. As a consequence, we obtain that ${\textbf{GL}\times\textbf{S5}} = \textbf{MGLB}$ and $\textbf{Grz}\times\textbf{S5} = \textbf{M}^+\textbf{GrzB}$.

The paper is structured as follows.
%We present our work in the following order. 
In \Cref{sec2}, we briefly recall
%provide a brief introduction to 
product logics and mm-logics. 
In \Cref{sec3}, we observe that $\textbf{GL}[n]\times\textbf{S5}$ is product matching, and that $\textbf{Grz}[n]\times\textbf{S5} = \textbf{M}^+\textbf{GrzB}[n]$ 
%is finitely axiomatizable 
for all $n\ge 1$. This we do by showing that the embedding $\textbf{Grz}[n]\hookrightarrow\textbf{GL}[n]$ lifts to an embedding $\textbf{M}^+\textbf{GrzB}[n]\hookrightarrow\textbf{MGLB}[n]$. %\color{blue} Maybe also mention that we lift the translation $\textbf{Grz}[n]\to\textbf{GL}[n]$. \color{black}
In \Cref{sec4}, we first outline our selective filtration technique from \cite{BM} in the absence of the monadic Barcan formula, and then provide necessary modifications to account for it. Finally, in \Cref{sec5}, we utilize the selective filtration technique from the previous section to prove the fmp for $\textbf{MGrzB}$, from which we derive the fmp for $\textbf{M}^+\textbf{GrzB}$ and $\textbf{MGLB}$. We use this to show that the embedding  $\textbf{Grz}\hookrightarrow\textbf{GL}$ extends to an embedding $\textbf{M}^+\textbf{GrzB}\hookrightarrow\textbf{MGLB}$. 
%\color{blue} Also mention that the embedding $\textbf{Grz}\hookrightarrow\textbf{GL}$ lifts to an embedding $\textbf{M}^+\textbf{GrzB}\hookrightarrow\textbf{MGLB}$. \color{red} 
The fmp of $\textbf{M}^+\textbf{GrzB}$ and $\textbf{MGLB}$ is further utilized to prove that $\textbf{GL}\times\textbf{S5}$ is product matching, and that $\textbf{Grz}\times\textbf{S5} = \textbf{M}^+\textbf{GrzB}$. \color{black} As a consequence, we obtain that $\textbf{GL}\times\textbf{S5}$ axiomatizes the one-variable fragment of $\textbf{QGLB}$, while  $\textbf{Grz}\times\textbf{S5}$ that of $\textbf{Q}^+\textbf{GrzB}$ (the latter logic is obtained by adding to $\textbf{QGrzB}$ the G\"odel translation of Casari's formula).
%, \color{red} where
%\begin{equation*}    \textbf{Q}^+\textbf{Grz}=\textbf{QGrz}+\Box\forall x(\Box(\Box p(x)\to \Box\forall xp(x))\to\Box\forall xp(x))\to\Box\forall xp(x).
%\end{equation*}
%\color{blue} (the latter system remains undefined). \color{black}

\section{Preliminaries}\label{sec2}

In this section, we briefly recall products of modal logics and their relationship to mm-logics. 
%Products of logics are studied extensively in the literature; for a comprehensive introduction, see \cite{mdim}.

%Recall the notion of product frames.

\begin{defn}\cite[p.~222]{mdim}
%\cite[p.~222, Sec.~5.1]{mdim}
\label{products}
    Let $\mathfrak{F}=(X_1,R_1)$ and $\mathfrak{G}=(X_2,R_2)$ be Kripke frames. The {\em product} of $\mathfrak{F}$ and $\mathfrak{G}$ is the frame $$\mathfrak{F}\times\mathfrak{G}=(X_1\times X_2, R_h, R_v)$$ where 
    %$X_1\times X_2$ is the Cartesian product of $X$ and $Y$, and 
    the relations $R_h, R_v$ are defined on $X_1\times X_2$ as follows for any $(x,u),(y,v)\in X_1\times X_2$:
    \begin{align*}
        (x,u)\mathrel{R_h}(y,v) &\iff x\mathrel{R_1}y\text{ and }u=v,\\
        (x,u)\mathrel{R_v}(y,v) &\iff x=y\text{ and }u\mathrel{R_2}v.
    \end{align*} 
\end{defn}

As we pointed out in the introduction, we will only work
%Given our interest in mm-logics, only products with $\textbf{S5}$ will be relevant to us. Hence, we work 
with product frames 
%of the form 
$\mathfrak{F}\times\mathfrak{G}$ where $\mathfrak{G}$ is an \textbf{S5}-frame (that is, $R_2$ is an equivalence relation). It is then straightforward to check that $R_v$ is an equivalence relation. 
%\color{red} (they don't mention anything about this, but it should be an easy verification) \color{blue} (add ref if this is in the yellow book). \color{black} 
We thus simplify our notation and denote $R_h$ by $R$ and $R_v$ by $E$.   
%; we denote the relation $R_h$ simply by $R$. 
%Henceforth, whenever products or product logics are mentioned, we invoke them in the context of \textbf{S5} as described here.
When we depict such frames, 
%In any of the diagrams shown in the paper, 
the equivalence relation $E$ is always 
%depicted 
drawn horizontally and the relation $R$ 
%is depicted 
vertically. 

%Given the semantic nature of the definition of product logics, it is natural to contemplate their axiomatization. In particular, investigating instances of \emph{product matching} is especially pertinent. We elaborate on this down below.

Let $\mathfrak{F}\times\mathfrak{G}=(X_1\times X_2,R,E)$ be a product frame. It is well known (see, e.g., \cite[p.~222]{mdim}) 
%\color{blue} [REF] \color{black} 
%By routine inspection, one can show 
that the following first-order properties hold in $\mathfrak{F}\times\mathfrak{G}$: 
%\color{red} done \color{blue} (instead of (1) and (2), should these be denoted by (LC) and (RC)?): \color{black} 
\begin{enumerate}[label={(LC)}]
    \item Left commutativity: $(\forall a,b,c\in X_1\times X_2)(a\mathrel{E}b \mbox{ and } b\mathrel{R} c) \implies (\exists d \in X_1 \times X_2)(a\mathrel{R} d \mbox{ and } d\mathrel{E} c)$. 
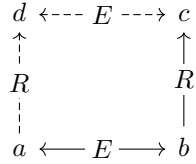
\label{left_comm}   
\[\begin{tikzcd}
	d && c \\
	\\
	a && b
	\arrow["E"{description}, tail reversed, dashed, from=1-1, to=1-3]
	\arrow["R"{description}, dashed, from=3-1, to=1-1]
	\arrow["E"{description}, tail reversed, from=3-1, to=3-3]
	\arrow["R"{description}, from=3-3, to=1-3]
\end{tikzcd}\]
\end{enumerate}

\begin{enumerate}[label={(RC)}]    
    \item Right commutativity: $(\forall a,b,c\in X_1\times X_2)(a\mathrel{R}b \mbox{ and } b\mathrel{E} c) \implies (\exists d \in X_1 \times X_2)(a\mathrel{E} d \mbox{ and } d\mathrel{R} c)$. 
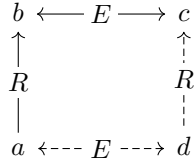
\label{right_comm}
\[\begin{tikzcd}
	b && c \\
	\\
	a && d
	\arrow["E"{description}, tail reversed, from=1-1, to=1-3]
	\arrow["R"{description}, from=3-1, to=1-1]
	\arrow["E"{description}, tail reversed, dashed, from=3-1, to=3-3]
	\arrow["R"{description}, dashed, from=3-3, to=1-3]
\end{tikzcd}\]
\end{enumerate}

To express these first-order properties modally, we follow the convention that a modal logic \textbf{L} is formulated in the propositional language $\mathcal{L}_\Diamond$ with the modality $\Diamond$, while \textbf{S5} in the propositional language $\mathcal{L}_\exists$ with the modality $\exists$. As usual, $\Box$ and $\forall$ abbreviate $\neg\Diamond\neg$ and $\neg\exists\neg$, respectively. Let $\mathcal{L}_{\Diamond\exists}$ denote the propositional language 
%containing the language of classical propositional logic and 
with two modalities $\Diamond$ and $\exists$. The {\em fusion} $\textbf{L}\otimes\textbf{S5}$ is the least set of formulae in  $\mathcal{L}_{\Diamond\exists}$ containing the theorems of $\textbf{L}$ and $\textbf{S5}$ and closed under the inference rules of Modus Ponens, Substitution, $\Box$-necessitation, and $\forall$-necessitation (see, e.g., \cite[p.~111]{mdim}).  %(note that $\Box$ and $\forall$ are abbreviations of $\neg\Diamond\neg$ and $\neg\exists\neg$, respectively). 
%\color{blue} (add a ref, also need to mention what $\Box$ and $\forall$ abbreviate). \color{black} 

%\color{blue} Need to mention what is meant by a 2-frame; add ref. 
%\color{red}
%Recall that 
The fusion $\textbf{L}\otimes\textbf{S5}$ is interpreted in a {\em $2$-frame} 
%is a structure of the form 
$\mathfrak{F}=(X,R,E)$, where 
%$X$ is a non-empty set, and 
$R$ is a binary relation and $E$ is an equivalence on $X$ (see, e.g., \cite[p.~21]{mdim}). In the $2$-frame $\mathfrak{F}$, the relation $R$ interprets the modality $\Diamond$, while the equivalence relation $E$ is reserved for the modality $\exists$.

Product frames are examples of $2$-frames, and it is straightforward to verify (see, e.g., \cite[p.~222]{mdim}) that left and right commutativity are definable in $\mathcal{L}_{\Diamond\exists}$ by 
%$\text{com}^l$ and $\text{com}^r$ formalize 
%are, in fact, modally expressible. 
%Let
\begin{align*}
    \text{com}^l&=\exists\Diamond p\to \Diamond\exists p,\\
    \text{com}^r&=\Diamond\exists p \to \exists\Diamond p.
\end{align*}

%For a Kripke complete propositional modal logic $\textbf{L}$, let $\textbf{L}\times\textbf{S5}$ denote the set of all formulas of $\mathcal{L}_{\Diamond\exists}$ valid in a product frame $\mathfrak{F}\times\mathfrak{G}$, where $\mathfrak{F}$ is a frame for $\textbf{L}$ and $\mathfrak{G}$ is a frame for $\textbf{S5}$. Following \color{blue} [REF], \color{black} we call $\textbf{L}\times\textbf{S5}$ a {\em product logic}.

\begin{defn}\cite[p.~126]{mdim}
    %For a Kripke complete propositional modal logic $\textbf{L}$, 
    %\begin{equation*}
    %    \textbf{L}\times\textbf{S5}=\text{Log}(\mathsf{C}_\textbf{L}),
    %\end{equation*}
    %where $\text{Log}(-)$ denotes the set of all valid $\mathcal{L}_{\Diamond\exists}$-formulae in the specified class. 
    %\color{blue} [REF] \color{black} 
    For a Kripke complete propositional modal logic $\textbf{L}$, let $\textbf{L}\times\textbf{S5}$ denote the set of all formulae of $\mathcal{L}_{\Diamond\exists}$ valid in all product frames $\mathfrak{F}\times\mathfrak{G}$, where $\mathfrak{F}$ is an $\textbf{L}$-frame and $\mathfrak{G}$ is an $\textbf{S5}$-frame. We call $\textbf{L}\times\textbf{S5}$ a {\em product logic}.
\end{defn}

Clearly $\text{com}^l$ and $\text{com}^r$ are provable in $\textbf{L}\times\textbf{S5}$, which motivates the following:

%\color{blue}
%Need to also define $\textbf{L}\times\textbf{S5}$. 
%\color{black}

\begin{defn}\cite[p.~223]{mdim}
%\cite[p.~223, Sec.~5.1]{mdim}
\label{matching}
    For a propositional modal logic \textbf{L}, let 
    \begin{equation*}
    [\textbf{L},\textbf{S5}] = (\textbf{L}\otimes \textbf{S5}) + \text{com}^l + \text{com}^r.
    \end{equation*}
    %where $\textbf{L}\otimes \textbf{S5}$ denotes the fusion of $\textbf{L}$ and $\textbf{S5}$. 
    We call $\textbf{L}\times\textbf{S5}$ {\em  product matching} provided $\textbf{L}\times \textbf{S5}=[\textbf{L},\textbf{S5}]$. 
\end{defn}

We next connect product logics with 
%ow formally define 
mm-logics. 

\begin{defn}\label{mm-logics} \cite[Def.~2.13]{GBm} 
    \begin{enumerate}[label={(\arabic*)}]
        \item For a propositional modal logic \textbf{L}, its {\em monadic extension} is the logic \label{monadic_extension}
            \begin{equation*}
                \textbf{ML} = (\textbf{L}\otimes\textbf{S5}) + \text{com}^l.
            \end{equation*}
        In particular, $\textbf{MK}$ is the monadic extension of the least  propositional normal modal logic $\textbf{K}$. 
        \item A {\em monadic modal logic} or {\em mm-logic}  is a set of formulae in the language $\mathcal{L}_{\Diamond\exists}$ containing \textbf{MK} and closed under the inference rules of Modus Ponens, Substitution, $\Box$-necessitation, and $\forall$-necessitation. 
    \end{enumerate}
\end{defn}

%\begin{rmk}
    In \cite[p.~432]{mdim}, the logic $\textbf{ML}$ is denoted by %$\textbf{L}\otimes\textbf{S5}+\text{com}^l$ is denoted by 
    $[\textbf{L},\textbf{S5}]^{\text{EX}}$. 
    %Thus, $\textbf{ML}=[\textbf{L},\textbf{S5}]^{\text{EX}}$.
%\end{rmk}
Indeed, many mm-logics can be realized as logics of \emph{expanding relativized products}, which can be seen as fragments of product logics (see, e.g., \cite[p.~432]{mdim}). 

%\begin{rmk}
    Recalling the {\em Barcan formula}  $$\Diamond\exists x P(x)\to\exists x \Diamond P(x)$$ and the {\em converse Barcan formula} $$\exists x\Diamond P(x)\to\Diamond\exists x P(x)$$ from predicate modal logic (see, e.g., \cite[pp.~244--245]{H&G}), 
    %\color{blue} (we probably need a different ref), \color{black} 
    we will refer to $\text{com}^r$ as the {\em monadic Barcan formula}, and to $\text{com}^l$ as the {\em converse monadic Barcan formula}.
%\end{rmk}
%Additionally, 
While each mm-logic contains the latter, it may not contain the former.
%monadic Barcan formula. 
%$\Diamond\exists p\to \exists\Diamond p$.
For an mm-logic \textbf{M}, we denote its extension with the monadic Barcan formula by $\textbf{MB}$. If \textbf{L} is a propositional modal logic, then its monadic extension with the monadic Barcan formula is the logic
%$\textbf{MLB}$ is given by
\begin{equation*}
    \textbf{MLB} = (\textbf{L}\otimes\textbf{S5}) + \text{com}^l + \text{com}^r.
\end{equation*}
Thus, $\textbf{MLB} = [\textbf{L},\textbf{S5}]$, and $(\textbf{L},\textbf{S5})$ is a product matching pair provided  $\textbf{L}\times\textbf{S5}=\textbf{MLB}$. 

In general, for a Kripke complete propositional modal logic $\textbf{L}$, we always have one inclusion $\textbf{MLB}\subseteq\textbf{L}\times\textbf{S5}$ (see, e.g., \cite[Prop.~3.8]{mdim}).
However, the other inclusion may not hold. When it does, we obtain
that $\textbf{MLB}$ axiomatizes the one-variable fragment of the predicate logic $\textbf{QLB}$ (see \cite[Thm.~13.8(1)]{gab98}, \cite[Thm.~3.21]{mdim} or, more generally, \cite[Thm.~5.4]{BM}). 
%While product matching is sufficient for $\textbf{MLB}$ to be the one-variable fragment of $\textbf{QLB}$, it is unclear if it is necessary.
%\color{blue} Mention that when it does hold, we get that $\textbf{MLB}$ axiomatizes the monadic fragment of $\textbf{QLB}$. In general, the connection to predicate modal logic should be discussed in some detail somewhere. \color{black}
%To prove the other inclusion, for any $\mathcal{L}_{\Diamond\exists}$-formula $\varphi$ not provable in $\textbf{MLB}$, we must find a product frame $\mathfrak{F}\times\mathfrak{G}$ for $\textbf{L}\times\textbf{S5}$ which refutes $\varphi$. 

%The frames for mm-logics are those 2-frames $\mathfrak{F}=(X,R,E)$ that satisfy \labelcref{left_comm}; additionally, $\mathfrak{F}$ validates the monadic Barcan formula iff $\mathfrak{F}$ satisfies \labelcref{right_comm} (see, e.g., \cite[p.~222]{mdim}) \color{red} probably better to rephrase this sentence since we have already mentioned this above\color{blue} (add ref for these two facts). \color{black} 

%\color{red}
We interpret mm-logics in $2$-frames $\mathfrak{F}=(X,R,E)$ that satisfy \labelcref{left_comm} (see, e.g., \cite[Def.~2.17]{GBm}), and refer to them as {\em $\textbf{MK}$-frames} (since they are exactly the $2$-frames validating $\textbf{MK}$). By the above discussion, $\mathfrak{F}$ validates the monadic Barcan formula iff $\mathfrak{F}$ satisfies \labelcref{right_comm}. 
%\color{red} We maintain the following convention while referring to mm-logics and their frames: if $\textbf{M}$ is any mm-logic and $\mathfrak{F}$ is a frame for $\textbf{M}$, then we call $\mathfrak{F}$ an $\textbf{M}$-frame. In particular, any frame for an mm-logic is an $\textbf{MK}$-frame. 
%\color{black}

%Given our goal to axiomatize $\textbf{Grz}\times\textbf{S5}$ and $\textbf{GL}\times\textbf{S5}$, as an intermediate step, we show that these product logics are indeed finitely axiomatizable in the presence of the finite depth formulae.  
%\color{blue}
%The ending might be a little too abrupt. 
%\color{black}

Since we are interested in the logics $\textbf{Grz}\times\textbf{S5}$ and $\textbf{GL}\times\textbf{S5}$, 
we recall the semantics for $\textbf{Grz}$ and $\textbf{GL}$. It is well known that 
\begin{itemize}
    \item $\textbf{Grz}$ is the logic of Noetherian posets, i.e., posets without infinite ascending chains (see, e.g., \cite[Cor.~5.52(i)]{CZ});
    \item $\textbf{GL}$ is the logic of Noetherian strict orders,  i.e., strict orders without infinite ascending chains (see, e.g., \cite[Cor.~5.47(i)]{CZ}).
\end{itemize}

%In the next section, 
We will also work with $\textbf{Grz}$ and $\textbf{GL}$ in the presence of the finite depth formulae, which are given by 
\begin{align*}
    bd_1&=\Diamond\Box p_1\to p_1,\\
    bd_n&=\Diamond(\Box p_n\wedge \neg bd_{n-1})\to p_n.
\end{align*}
These formulae express the finiteness of the depth of a transitive frame. To this end, we recall: 
%definition.

\begin{defn}\label{depth_defn}
    Let $\mathfrak{F}=(X,R,E)$ be an \textbf{MK}-frame with $R$ transitive. 
    \begin{enumerate}[label={(\arabic*)}]
        \item A chain $x_1\mathrel{R}\dots\mathrel{R} x_n$ in $\mathfrak{F}$ is a {\em strict chain} if $x_{i+1}\mathrel{\cancel{R}}x_i$ for any $i\in\{1,\dots,n-1\}$. The {\em length} of a strict chain refers to the cardinality of its underlying set. \label{strict_chain}
        
        \item The frame $\mathfrak{F}$ is of {\em depth $n$}, written $d(\mathfrak{F})=n$, provided it contains a strict chain of length $n$ and no other strict chain in $\mathfrak{F}$ has greater length. 
        %If $\mathfrak{F}$ is not of depth $n$ for any $n<\omega$, then it has {\em depth $\omega$} \color{red} we don't use it anywhere---should I delete it? \color{blue} (do we ever use this?). \color{black} We denote the depth of $\mathfrak{F}$ by $d(\mathfrak{F})$. 
        \label{depth_n}
    \end{enumerate}
\end{defn}

The following result is a direct adaptation of a well-known result to the monadic setting: 

\begin{prop} (see, e.g., \cite[Prop.~3.44]{CZ})
    Let $\mathfrak{F}=(X,R,E)$ be an \textbf{MK}-frame with $R$ transitive. Then $\mathfrak{F}\models bd_n$ iff $d(\mathfrak{F})\leq n$. 
\end{prop}

For any logic \textbf{L}, we denote by $\textbf{L}[n]$ the extension of \textbf{L} with the finite depth formula $bd_n$. Hence, $\textbf{Grz}[n]$ is the logic of those $\textbf{Grz}$-frames that have depth at most $n$, and $\textbf{GL}[n]$ is defined similarly. 
%the logic of those $\textbf{GL}$-frames that have depth at most $n$. 
In the next section we will show that  $\textbf{GL}[n]\times\textbf{S5}$ and $\textbf{Grz}[n]\times\textbf{S5}$ are finitely axiomatizable for each $n\geq 1$.

\section{Axiomatizing \texorpdfstring{$\textbf{Grz}[n]\times\textbf{S5}$}{Grz[n] x S5}}\label{sec3}

%\color{blue} I wonder whether the definitions of \textbf{Grz}, \textbf{GL}, and their finite depth versions should already be given in the intro. Then this section can start where it started earlier. 
%Need to transition to $\bf Grz$ and $\bf GL$ and their semantics. It might be better to move the business about ${\bf Grz}[n]$ and ${\bf GL}[n]$ to the next section. 
%\color{black} 

%\color{blue}
%Probably best to start the section with the next paragraph. 
%\color{black}

As we pointed out in the introduction, product matching is not a common phenomenon. For our purposes, we recall that neither $\textbf{Grz}\times\textbf{S5}$ nor $\textbf{Grz}[n]\times\textbf{S5}$ ($n\geq 2$) is product matching (see \cite[Thm.~5.17]{mdim}). We will show that these product logics are finitely axiomatizable. In fact, we will prove that they can be axiomatized by the G\"odel translation of the monadic Casari formula $$\forall((p\to\forall p)\to\forall p)\to\forall p,$$ which plays an important role in lifting the well-known translations 
\[
\textbf{IPC} \hookrightarrow \textbf{Grz} \hookrightarrow \textbf{GL}
\]
to the monadic setting \cite{GBm}. In this section, we concentrate on $\textbf{Grz}[n]\times\textbf{S5}$ and in the next on $\textbf{Grz}\times\textbf{S5}$. 

Recall from \crefdefpart{mm-logics}{monadic_extension} that $\textbf{MGrz}$ is the monadic extension of $\textbf{Grz}$, and that adding the monadic Barcan formula to $\textbf{MGrz}$ results in the mm-logic $\textbf{MGrzB}$. Following \cite[Def.~4.5]{GBm}, we set: 

\begin{defn} Let 
\begin{align*}
    \textbf{M}^+\textbf{Grz}&=\textbf{MGrz}+\Box\forall(\Box(\Box p\to \Box\forall p)\to\Box\forall p)\to\Box\forall p,\\
    \textbf{M}^+\textbf{GrzB}&=\textbf{MGrzB}+\Box\forall(\Box(\Box p\to \Box\forall p)\to\Box\forall p)\to\Box\forall p.
\end{align*}    
\end{defn}
\begin{rmk}
    The defining formula above is equivalent to the G\"{o}del translation of the monadic Casari formula \cite[p.~441]{GBm}. For simplicity, we will refer to it as the {\em Casari formula}.
\end{rmk}

%We recall the G\"{o}del translation of the monadic Casari formula, $\Box\forall(\Box(\Box p\to \Box\forall p)\to\Box\forall p)\to\Box\forall p$. We will refer to this simply as the Casari formula. Define 

%Our primary aim will be to show that $\textbf{Grz}\times\textbf{S5}=\textbf{M}^+\textbf{GrzB}$. However, we find it instructive to prove that $\textbf{M}^+\textbf{GrzB}[n]=\textbf{Grz}[n]\times\textbf{S5}$, as this will provide an overview of our strategy.  

%The Casari formula was extensively studied in \cite{GBm}, where it was given a semantic characterization in the setting of descriptive frames. In what follows, we present a similar characterization that clarifies the content of the formula.

To give a semantic characterization of the Casari formula, we recall that for a binary relation $S$ on a set $X$, the $S$-image of $x \in X$ is the set $S[x]:=\{y\in X\mid x\mathrel{S}y\}$. If $S$ is an equivalence relation, then we refer to $S[x]$ as an {\em $S$-cluster}.

\begin{defn}\cite[Def.~3.6(3)]{GBm} \label{clean}
    Let $\mathfrak{F}=(X,R,E)$ be an \textbf{MK}-frame. For $x \in X$, we say that the $E$-cluster $E[x]$ is {\em clean} provided $u,v\in E[x]$ and $uRv$ imply $u = v$. Otherwise, we say that the cluster $E[x]$ is {\em dirty}. 
    %\begin{equation*}
    %    x\mathrel{E}y, x\mathrel{R}y\implies x=y.
    %\end{equation*}
\end{defn}

\begin{exm} \label{casari_rmk}
A typical example of an \textbf{MGrz}-frame with a dirty cluster is the following frame:  \label{casari_rmk1}

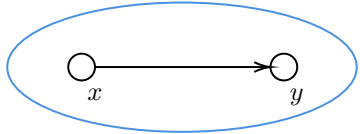
    \begin{center}
    \tikzset{every picture/.style={line width=0.75pt}} %set default line width to 0.75pt        

    \begin{tikzpicture}     [x=0.75pt,y=0.75pt,yscale=-1,xscale=1, scale=0.75]

    %Shape: Ellipse [id:dp021421804850109027] 
    \draw  [color={rgb, 255:red, 74; green, 144; blue, 226 }  ,draw opacity=1 ] (41,148.1) .. controls (41,124.08) and (93.61,104.6) .. (158.5,104.6) .. controls (223.39,104.6) and (276,124.08) .. (276,148.1) .. controls (276,172.12) and (223.39,191.6) .. (158.5,191.6) .. controls (93.61,191.6) and (41,172.12) .. (41,148.1) -- cycle ;
    %Shape: Circle [id:dp20738343997713837] 
    \draw   (82,148) .. controls (82,143.03) and (86.03,139) .. (91,139) .. controls (95.97,139) and (100,143.03) .. (100,148) .. controls (100,152.97) and (95.97,157) .. (91,157) .. controls (86.03,157) and (82,152.97) .. (82,148) -- cycle ;
    %Shape: Circle [id:dp4766861728645406] 
    \draw   (218,148) .. controls (218,143.03) and (222.03,139) .. (227,139) .. controls (231.97,139) and (236,143.03) .. (236,148) .. controls (236,152.97) and (231.97,157) .. (227,157) .. controls (222.03,157) and (218,152.97) .. (218,148) -- cycle ;
    %Straight Lines [id:da5328173386993165] 
    \draw    (100,148) -- (216,148) ;
    \draw [shift={(218,148)}, rotate = 180] [color={rgb, 255:red, 0; green, 0; blue, 0 }  ][line width=0.75]    (10.93,-3.29) .. controls (6.95,-1.4) and (3.31,-0.3) .. (0,0) .. controls (3.31,0.3) and (6.95,1.4) .. (10.93,3.29)   ;

    % Text Node
    \draw (93,160.4) node [anchor=north west][inner sep=0.75pt]    {$x$};
    % Text Node
    \draw (229,160.4) node [anchor=north west][inner sep=0.75pt]    {$y$};

    \end{tikzpicture}
    \end{center}

    Here, $X=\{x,y\}$, $R=\{(x,x),(x,y),(y,y)\}$, and $E=X^2$ (where the circles represent reflexive points, the black arrow indicates the non-reflexive relation $(x,y)$, and the blue curve depicts the equivalence relation $E$). %\color{blue} Can you please make the diagram a little smaller and also explain what the blue circle represents? \color{black} 
    Observe that this frame is in fact an $\textbf{MGrzB}$-frame.  
\end{exm}

The next result was proved in \cite[Lem.~4.8]{GBm} for descriptive $\textbf{MGrz}$-frames, but since $\textbf{Grz}$-frames have no infinite ascending $R$-chains, the same proof yields the following:

\begin{thm}\label{casari}
    Let $\mathfrak{F}=(X,R,E)$ be an $\textbf{MGrz}$-frame. Then $\mathfrak{F}$ validates the Casari formula %$\mathfrak{F}\models\Box\forall(\Box(\Box p\to \Box\forall p)\to\Box\forall p)\to\Box\forall p$ 
    iff $\mathfrak{F}$ has clean $E$-clusters.
\end{thm}
%\begin{proof}
%    See \cite[Lem.~4.8]{GBm}. 
%\end{proof}

\begin{rmk} \label{casari_rmk2}
Since each $(\textbf{Grz}\times\textbf{S5})$-frame has clean clusters, \Cref{casari} shows that each such frame validates the Casari formula. On the other hand, the frame depicted in \cref{casari_rmk} is an $\textbf{MGrzB}$-frame that refutes the Casari formula. Since the depth of this frame is $2$, we obtain yet another proof that neither $\textbf{Grz}\times\textbf{S5}$ nor $\textbf{Grz}[n]\times\textbf{S5}$, for $n\geq 2$, is product matching.
%a frame for $\textbf{Grz}[n]\times\textbf{S5}$ is a frame for $\textbf{M}^+\textbf{GrzB}[n]$. For $n\geq 2$, it is easy to see that $\textbf{M}^+\textbf{GrzB}[n]$ is a proper extension of $\textbf{MGrzB}[n]$ (the example above is a frame for \textbf{MGrzB} that does not have clean $E$-clusters). Therefore, \Cref{casari} provides another proof that $\textbf{Grz}[n]\times\textbf{S5}$ is not product matching whenever $n\geq 2$. The same observations are also applicable to $\textbf{Grz}\times\textbf{S5}$.
\end{rmk}

%Our approach will require the completeness of $\textbf{M}^+\textbf{GrzB}$; this will follow from the fact that it has the fmp, which we prove in the next section. Nonetheless, we have the fmp (and hence completeness) of $\textbf{M}^+\textbf{GrzB}[n]$ due to the following. 

To prove that $\textbf{Grz}[n]\times\textbf{S5} = \textbf{M}^+\textbf{GrzB}[n]$, we first show that the well-known embedding $\textbf{Grz}[n] \hookrightarrow\textbf{GL}[n]$ extends to yield the embeddings $\textbf{M}^+\textbf{Grz}[n] \hookrightarrow\textbf{MGL}[n]$ and $\textbf{M}^+\textbf{GrzB}[n] \hookrightarrow\textbf{MGLB}[n]$, and then utilize that $(\textbf{GL}[n],\textbf{S5})$ is a product matching pair. To see the latter, we recall the notion of Horn axiomatizable unimodal logics. 

\begin{defn}\cite[p.~228]{mdim}\label{Horn_axiom}
    A Kripke complete unimodal logic $\textbf{L}$ is \emph{Horn axiomatizable} if the class of $\textbf{L}$-frames can be axiomatized by closed formulae or by the formulae that correspond to first-order sentences of the form
    \begin{equation*}
        \forall x\forall y\forall \Vec{z} \big(\psi(x,y,\Vec{z})\to R(x,y)\big),
    \end{equation*}
where $\psi(x,y,\Vec{z})$ is a positive formula.
\end{defn}

The following 
%is the most general 
result justifies the importance of Horn axiomatizable unimodal logics.
%known about product matching with \textbf{S5}.

\begin{thm}\label{Horn}\cite[Thm.~5.9]{mdim}
    If $\textbf{L}_1$ and $\textbf{L}_2$ are Horn axiomatizable unimodal logics, then $\textbf{L}_1\times\textbf{L}_2$ is product matching. In particular, if $\textbf{L}$ is Horn axiomatizable, then $\textbf{L}\times\textbf{S5}$ is product matching.
\end{thm}

\begin{rmk}
%    \item[]
%    \item 
    \Cref{Horn} implies that $\textbf{K}\times\textbf{S5}$, $\textbf{K4}\times\textbf{S5}$, $\textbf{S4}\times\textbf{S5}$, and $\textbf{S5}\times\textbf{S5}$ are all product matching. However, the same conclusion cannot be made about $\textbf{Grz}\times\textbf{S5}$ and $\textbf{GL}\times\textbf{S5}$ because neither $\textbf{Grz}$ nor $\textbf{GL}$ is Horn axiomatizable (in fact, their classes of frames are not even first-order definable; see, e.g., \cite[Thm.~6.7, 6.8]{CZ}). Indeed, we already saw earlier in the section that $\textbf{Grz}\times\textbf{S5}$ is not product matching. Nevertheless, as we will see in \crefdefpart{main}{product_matching_GL},
    %\color{blue} (add exact forward reference), \color{black} 
    $\textbf{GL}\times\textbf{S5}$ is product matching, but the proof requires a different approach.
    %
    %\item 
    %Note that due to the compactness of first order logic, frames for $\textbf{GL}$ and $\textbf{Grz}$ are not first order axiomatizable (\cite[Thm.~6.7, Thm.~6.8]{CZ}); so, $\textbf{GL}$ and $\textbf{Grz}$ are not Horn axiomatizable. Therefore, \Cref{Horn} is not applicable directly in the case of these logics. 
    % 
    %\item 
%\end{enumerate}
\end{rmk}

Since $\textbf{Grz}[n]\times\textbf{S5}$ is not product matching for $n \ge 2$, we conclude from \Cref{Horn} that $\textbf{Grz}[n]$ is not Horn axiomatizable. %This is, in part, due to the fact that Horn formulae are preserved under products, while the depth of partial orders is not. However, we know that 
On the other hand, each $\textbf{GL}[n]$ is Horn axiomatizable. To see this, note that $\textbf{GL}[n]=\textbf{K4}+\neg\Diamond^n \top$ for each $n \ge 1$. Indeed, for a Kripke frame $\mathfrak F=(X,R)$, we have $\mathfrak F\models\neg\Diamond^n\top$ iff $R^n=\varnothing$ (see, e.g., \cite[Prop. 3.33]{venema}). In particular, $\mathfrak F\models\neg\Diamond^n\top$ implies that every point in $\mathfrak F$ is irreflexive.
%: for if $x\mathrel{R}x$, then $x\mathrel{R^n}x$. For 
When $\mathfrak F$ is transitive, the condition $R^n=\varnothing$ 
%it is straightforward to verify 
also implies that 
%$(X,R)\models\neg\Diamond^n\top$ implies 
$d(\mathfrak{F})\leq n$, yielding that $\mathfrak F \models \textbf{K4}+\neg\Diamond^n\top$ iff $\mathfrak F$ is a transitive irreflexive frame of depth $\le n$. Since $\textbf{GL}[n]$ is also the logic of such frames (see, e.g., \cite[p.~429]{kracht}), $\textbf{GL}[n]=\textbf{K4}+\neg\Diamond^n \top$. 
%, i.e., $(X,R)$ has depth at most $n$. 
Thus, $\textbf{GL}[n]$ is Horn axiomatizable because so is 
%$\textbf{K4}+\neg\Diamond^n \top$ 
$\textbf{K4}+\neg\Diamond^n\top$. 
%Recalling that $\textbf{GL}[n]$ is Kripke complete has the same frames as $\textbf{K4}+\neg\Diamond^n\top$ 
%logics whose frames are transitive irreflexive Kripke frames of depth $\le n$ 
%\color{blue} [REF], \color{black} 
%is the logic of all transitive frames (see, e.g., \cite[Thm.~5.16]{CZ}), the logic $\textbf{K4}+\neg\Diamond^n\top$ is canonical, and hence complete with respect to irreflexive and transitive frames with depth at most $n$. However, $\textbf{GL}[n]$ is complete and also has the same class of frames.\color{black} Thus, we have the following:
%we conclude that $\textbf{GL}[n]=\textbf{K4}+\neg\Diamond^n \top$ for each $n \ge 1$. 
This together with \Cref{Horn} yields:

%\begin{thm}\label{GL_horn}
    %Let $n<\omega$ and 
%    $\textbf{GL}[n]=\textbf{K4}+\neg\Diamond^n \top$ for each $n \ge 1$. 
%    \label{Abashidze}
    %Then $\textbf{GL}[n]=\textbf{L}_n$. 
%\end{thm}

%Putting \Cref{Horn,GL_horn} 
%\color{blue} (why is ``Theorems'' added manually?) \color{black} 
%together, we obtain: 
%\color{blue} (the code needs adjusting so theorem is not mentioned twice): \color{black}

\begin{cor}\label{matching_GLn}
    %For any $n<\omega$, the logic $\textbf{GL}[n]$ is Horn axiomatizable and hence 
    $\textbf{GL}[n]\times\textbf{S5}$ is product matching for each $n\ge 1$.
\end{cor}

%\begin{proof}
%    Since $\textbf{L}_n$ is Horn axiomatizable, a direct application of \Cref{Horn} completes the proof.
%\end{proof} 

To show that the translation $\textbf{Grz}[n] \hookrightarrow\textbf{GL}[n]$ lifts to the translation $\textbf{M}^+\textbf{Grz}[n] \hookrightarrow\textbf{MGL}[n]$, we extend the reflexivization and irreflexivization operations from Kripke frames to \textbf{MK}-frames.

%It is well known that the key difference between the frames for \textbf{Grz} and \textbf{GL} is reflexivity \cite[p.~98, Thm.~3.88]{CZ}. This semantic connection has an analog in the monadic setting, which will prove instrumental for our purposes. We elaborate on this below.

Let $\mathfrak{F}=(X,R)$ be a Kripke frame. We recall (see, e.g., \cite[p.~98]{CZ})  
%\color{blue} [REF] \color{black} 
that the {\em reflexivization} of $\mathfrak F$ is the Kripke frame  $\mathfrak{F}^r=(X,R^r)$, where $R^r$ is the reflexive closure of $R$; i.e., 
\[
x\mathrel{R^r}y\iff x\mathrel{R}y \text{ or } x=y.
\]
We also recall 
%\cite[p.~98]{CZ} \color{blue} [REF] \color{black} 
that the {\em irreflexivization} of $\mathfrak F$ is the Kripke frame $\mathfrak{F}^{ir}=(X,R^{ir})$, where $R^{ir}$ is the irreflexive fragment of $R$; i.e., 
\[
x\mathrel{R^{ir}}y\iff x\mathrel{R}y \text{ and } x\neq y.
\]
These notions naturally extend to \textbf{MK}-frames:

\begin{defn}
    Let $\mathfrak{F}=(X,R,E)$ be an \textbf{MK}-frame. The {\em reflexivization} of $\mathfrak F$ is the frame $\mathfrak{F}^r=(X,R^r,E)$ and the {\em irreflexivization} of $\mathfrak F$ is the frame $\mathfrak{F}^{ir}=(X,R^{ir},E)$. 
\end{defn}

%Similarly, if $\mathfrak{F}=(X,R,E)$ is any 2-frame, we let $\mathfrak{F}^r=(X,R^r,E)$.

\begin{prop}\label{refl_irrefl}
    Let $\mathfrak{F}=(X,R,E)$ be an \textbf{MK}-frame.  
    \begin{enumerate}[label={(\arabic*)}]
        \item $\mathfrak{F}^r$ is also an \textbf{MK}-frame. Moreover, if $\mathfrak{F}$ validates $\textbf{B}$, then so does $\mathfrak{F}^r$. \label{refl}
        \item If $\mathfrak{F}$ has clean clusters, then $\mathfrak{F}^{ir}$ is also an \textbf{MK}-frame (with clean clusters). Moreover, if $\mathfrak{F}$ validates \textbf{B}, then so does  $\mathfrak{F}^{ir}$. \label{irrefl}
    \end{enumerate}

\end{prop}
\begin{proof}
     We need to check that the necessary commutativity conditions are satisfied. We will only verify left commutativity for the two claims above. The right commutativity is verified similarly.  
     
    \labelcref{refl} Suppose $x\mathrel{E}y$ and $y\mathrel{R^r}z$ for some $x,y,z\in X$. 
    We must find $u\in X$ such that $x\mathrel{R^r}u$ and $u\mathrel{E}z$. 
    From $y\mathrel{R^r}z$ it follows that $y=z$ or $y\mathrel{R}z$. If $y=z$, then we take $u=x$; and if $y\mathrel{R}z$, such a $u$ exists due to the fact that $\mathfrak{F}$ satisfies left commutativity. 

    \labelcref{irrefl} Suppose $x\mathrel{E}y$ and $y\mathrel{R^{ir}}z$ for some $x,y,z\in X$. We must find $u\in X$ such that $x\mathrel{R^{ir}}u$ and $u\mathrel{E}z$. From $y\mathrel{R^{ir}}z$ it follows that  $y\mathrel{R}z$ and $y\neq z$. Therefore, $x\mathrel{E}y$ and $y\mathrel{R}z$, and since $\mathfrak{F}$ satisfies left commutativity, there is $u\in X$ such that $x\mathrel{R}u$ and $u\mathrel{E}z$. Because  $\mathfrak{F}$ has clean clusters, $y\mathrel{R}z$ and $y\neq z$ imply that $y\mathrel{\cancel{E}}z$, and so  
    $x\mathrel{\cancel{E}}z$. Thus, $x\neq u$, and hence $x\mathrel{R^{ir}}u$.
    %
    %\color{blue} Add short proof for both items and explain why in (2) clean clusters are needed. \color{black}  
    %The relation $E$ is an equivalence. Therefore, it suffices to show that $\mathfrak{F}^r$ satisfies the required commutativity condition(s). This is easily checked. 
\end{proof}

\begin{rmk}
The assumption in \crefdefpart{refl_irrefl}{irrefl} 
%\color{blue} (fix the label) \color{black}
    is necessary. Indeed, 
    %that $\mathfrak{F}$ has clean clusters, as otherwise the claim in \Cref{irrefl} is false. For example, consider 
in the \textbf{MK}-frame shown in \cref{casari_rmk},
%\crefdefpart{casari_rmk}{casari_rmk1}. In this case, $R^{ir}=\{(x,y)\}$, and so 
$x\mathrel{E}y$ and $x\mathrel{R^{ir}}y$, but $R^{ir}[y]=\varnothing$. Hence, $\mathfrak{F}^{ir}$ does not satisfy left commutativity.
%
%    \item 
%    If $\mathfrak{F}\times\mathfrak{G}$ is a $\textbf{K}\times\textbf{S5}$-frame, then $(\mathfrak{F}\times\mathfrak{G})^{ir}=\mathfrak{F}^{ir}\times\mathfrak{G}$. However, we will not need this fact in our proofs.
%\end{enumerate}
\end{rmk}

%Let $\mathfrak{F}=(X,R)$ be a Kripke frame. Then we let $\mathfrak{F}^{ir}=(X,R^{ir})$, where $R^{ir}$ is the irreflexivization of $R$, i.e., $x\mathrel{R^{ir}}y\iff x\neq y \text{ and } x\mathrel{R}y$. Similarly, if $\mathfrak{F}=(X,R,E)$ is a 2-frame, then we take $\mathfrak{F}^{ir}=(X,R^{ir},E)$.

%\begin{prop}\label{irrefl}
%    Let $\mathfrak{F}=(X,R,E)$ be an \textbf{MK}-frame with clean clusters. Then $\mathfrak{F}^{ir}=(X,R^{ir},E)$ is also an \textbf{MK}-frame with clean clusters. Similarly, if $\mathfrak{F}$ is an \textbf{MKB}-frame with clean clusters, then so is $\mathfrak{F}^{ir}$.
%\end{prop}
%\begin{proof}
%    As in \Cref{refl}, we only need to verify the necessary commutativity condition(s).
%\end{proof}

%Recall the splitting translation $(-)^+\colon\textbf{Grz}\to\textbf{GL}$, 
As is customary, for a formula $\varphi$, we write $\Box^+\varphi := \varphi\wedge\Box\varphi$. The same proof as in \cite[Lem.~3.86]{CZ} yields:
%obtained by replacing any instance of $\Box$ in a formula with $\Box^+$, where $\Box^+\varphi=\varphi\wedge\Box\varphi$. 

%\begin{lem}\cite[p.~98, Lem.~3.86]{CZ}
%    Let $\mathfrak{F}=(X,R,E)$ be an \textbf{MK}-frame and $v$ be any valuation on $\mathfrak{F}$. Let $\mathfrak{M}=(\mathfrak{F},v)$ be the corresponding model. Then for any $x\in X$,
%    \begin{equation*}
%        \mathfrak{M},x\models_v\varphi^+\iff\mathfrak{M}^r,x\models_v\varphi.
%    \end{equation*}
%\end{lem}

\begin{lem}\label{refl_valid}
    Let $\mathfrak{F}=(X,R,E)$ be an \textbf{MK}-frame. Then for any $\mathcal{L}_{\Diamond\exists}$-formula $\varphi$,  %\color{blue} (say what $\varphi$ is) \color{black}
    \begin{equation*}    \mathfrak{F}\models\varphi^+\iff\mathfrak{F}^r\models\varphi.
    \end{equation*}
\end{lem}

%We note the semantic connection between $\textbf{MGLB}[n]$ and $\textbf{M}^+\textbf{GrzB}[n]$. 
The next result is a  generalization of \cite[Lem.~4.10]{GBm}. %\color{blue} Barcan plays no role in the Kristina paper, and the next lemma has an obvious analog without Barcan. So things need to be explained better. \color{black}

\begin{lem}\label{frames}
    \begin{enumerate}[label={(\arabic*)}]
        \item[]
        \item If $\mathfrak{F}$ is an $\textbf{MGL}[n]$-frame, then $\mathfrak{F}^r$ is an $\textbf{M}^+\textbf{Grz}[n]$-frame. \label{frames3}
        
        \item If $\mathfrak{F}$ is an $\textbf{MGLB}[n]$-frame, then $\mathfrak{F}^r$ is an $\textbf{M}^+\textbf{GrzB}[n]$-frame. \label{frames1}

        \item If $\mathfrak{F}$ is an $\textbf{M}^+\textbf{Grz}[n]$-frame, then $\mathfrak{F}^{ir}$ is an $\textbf{MGL}[n]$-frame. \label{frames4}

        \item If $\mathfrak{F}$ is an $\textbf{M}^+\textbf{GrzB}[n]$-frame, then $\mathfrak{F}^{ir}$ is an $\textbf{MGLB}[n]$-frame. \label{frames2}
    \end{enumerate}
\end{lem}
\begin{proof}
    We only prove \labelcref{frames1} and \labelcref{frames2}. The proofs for \labelcref{frames3} and \labelcref{frames4} are analogous. 
    
    \labelcref{frames1} Suppose $\mathfrak{F}=(X,R,E)$ is an $\textbf{MGLB}[n]$-frame. By \crefdefpart{refl_irrefl}{refl}, $\mathfrak{F}^r$ is an \textbf{MKB}-frame. Moreover, $R^r$ is a partial order since $R$ is a strict order. Furthermore,  $\mathfrak{F}^r$ has clean clusters since $\mathfrak{F}$ does (as $\mathfrak{F}$ is an $\textbf{MGL}$-frame). We show that $d(\mathfrak{F}^r)\leq n$. Suppose $x_1\mathrel{R^r}\dots\mathrel{R^r} x_k$ is a strict chain (see \crefdefpart{depth_defn}{strict_chain}). Then $x_{i+1}\mathrel{\cancel{R^r}}x_i$, and so  $x_{i+1}\mathrel{\cancel{R}}x_i$ for each $i\in\{1,.\dots,k-1\}$. 
    %since $R\subseteq R^r$. 
    Hence, $x_1\mathrel{R}\dots\mathrel{R}x_k$ is a strict chain in $(X,R)$. Since $d(\mathfrak{F})\leq n$, we obtain $k\leq n$, which yields 
    that $d(\mathfrak{F}^r)\leq n$. Thus, $\mathfrak{F}^r$ is a $\textbf{M}^+\textbf{GrzB}[n]$-frame. 

    \labelcref{frames2} Suppose $\mathfrak{F}=(X,R,E)$ is an $\textbf{M}^+\textbf{GrzB}[n]$-frame. By \Cref{casari}, $\mathfrak{F}$ has clean clusters, so by \crefdefpart{refl_irrefl}{irrefl}, 
    %\color{blue} (fix the label), \color{black} 
    $\mathfrak{F}^{ir}$ is an \textbf{MKB}-frame (with clean clusters). Moreover, $R^{ir}$ is a strict order since $R$ is a partial order. We show that  $d(\mathfrak{F}^{ir})\leq n$. Suppose $x_1\mathrel{R^{ir}}\dots\mathrel{R^{ir}} x_k$ is a strict chain. Then $x_1\mathrel{R}\dots\mathrel{R} x_k$. To see that this is a strict chain, let $x_{i+1}\mathrel{R}x_i$ for some $i$. 
    %Since $x_{i}\mathrel{R^{ir}} x_{i+1}\implies x_i\neq x_{i+1}$, we must have $x_i\mathrel{R}x_{i+1}\mathrel{R}x_i$, so  
    Then $x_i=x_{i+1}$ by antisymmetry of $R$, which is a contradiction since $x_{i}\mathrel{R^{ir}} x_{i+1}\implies x_i\neq x_{i+1}$. Therefore,  $d(\mathfrak{F})\leq n$ implies that $k\leq n$. Thus, $d(\mathfrak{F}^{ir})\leq n$, and hence $\mathfrak{F}^{ir}$ is an $\textbf{MGLB}[n]$-frame. 
\end{proof}

 By \cite[Thm.~4.12]{GBm}, $\varphi\in\textbf{M}^+\textbf{Grz}\iff\varphi^+\in\textbf{MGL}$. To obtain a similar result for the logics in the previous lemma, we require the following:
 %extend this further, we demonstrate that the translation $(-)^+\colon\textbf{Grz}\to\textbf{GL}$ lifts to the logics $\textbf{M}^+\textbf{GrzB}[n]$ and $\textbf{MGLB}[n]$. 

\begin{thm}\cite[Sec.~4.10, 4.11]{locally_finite}\label{locally_finite}
    The logics $\textbf{MGrz}[n]$ and $\textbf{MGL}[n]$ are locally tabular for each $n\ge 1$.
    %$n<\omega$.
\end{thm}

\begin{rmk}\label{p-morphicimage} 
     We recall that a $2$-frame $\mathfrak F' = (X',R',E')$ is a {\em p-morphic image} of a $2$-frame $\mathfrak F = (X,R,E)$ provided there is an onto map $f\colon X \to X'$ such that $f(R[x])=R'[f(x)]$ and $f(E[x])=E'[f(x)]$ for each $x\in X$ (see, e.g., \cite[p.~22-23]{mdim}). By \Cref{locally_finite}, for any $\mathcal{L}_{\Diamond\exists}$-formula $\varphi$, if $\mathfrak{F}$ is an $\textbf{MGrzB}[n]$-frame refuting $\varphi$, then there is a finite $\textbf{MGrzB}[n]$-frame $\mathfrak{G}$ that is a p-morphic image of $\mathfrak{F}$ and also refutes $\varphi$. Similar observations can be made for $\textbf{M}^+\textbf{GrzB}[n]$ and $\textbf{MGLB}[n]$. This
    %e following observation 
    will be used %given its relevance 
    in \Cref{MGrzB_fmp,fmp}.
\end{rmk}

\begin{thm}\label{translation} %\color{blue} Make this similar to \cref{frames}. \color{black} 
For an $\mathcal{L}_{\Diamond\exists}$-formula $\varphi$,
\begin{enumerate}[label={(\arabic*)}]
    \item $\varphi\in\textbf{M}^+\textbf{Grz}
    [n]\implies \varphi^+\in\textbf{MGL}[n]$. \label{translation1}

    \item $\varphi\in\textbf{M}^+\textbf{GrzB}[n]\implies \varphi^+\in\textbf{MGLB}[n]$. \label{translation2}

     \item $\varphi^+\in\textbf{MGL}[n]\implies \varphi\in\textbf{M}^+\textbf{Grz}[n]$. \label{translation3}
    
    \item $\varphi^+\in\textbf{MGLB}[n]\implies \varphi\in\textbf{M}^+\textbf{GrzB}[n]$. \label{translation4}

\end{enumerate}
\end{thm}
\begin{proof}
    Once again, we only prove \labelcref{translation2} and \labelcref{translation4}. The proofs for \labelcref{translation1} and \labelcref{translation3} are entirely similar.
    
    \labelcref{translation2} Let $\varphi^+\notin\textbf{MGLB}[n]$. By \Cref{locally_finite}, there is a finite $\textbf{MGLB}[n]$-frame $\mathfrak{F}$ such that $\mathfrak{F}\not\models\varphi^+$. By \Cref{refl_valid}, $\mathfrak{F}^r\not\models\varphi$. By \crefdefpart{frames}{frames1}, $\mathfrak{F}^r$ is an $\textbf{M}^+\textbf{GrzB}[n]$-frame. 
    %From $\mathfrak{F}^r\not\models\varphi$, we deduce that 
    Thus, $\varphi\notin\textbf{M}^+\textbf{GrzB}[n]$. 
    
    \labelcref{translation4} Let  $\varphi\notin\textbf{M}^+\textbf{GrzB}[n]$. By \Cref{locally_finite}, there is a finite $\textbf{M}^+\textbf{GrzB}[n]$-frame $\mathfrak{F}$ such that $\mathfrak{F}\not\models\varphi$. By \crefdefpart{frames}{frames2}, $\mathfrak{F}^{ir}$ is an $\textbf{MGLB}[n]$-frame, and it is easy to see that  $(\mathfrak{F}^{ir})^r=\mathfrak{F}$. Therefore, by \Cref{refl_valid},  $\mathfrak{F}^{ir}\not\models\varphi^+$. Thus, $\varphi^+\notin\textbf{MGLB}[n]$.
\end{proof}

%Now, we claim that $\textbf{MGLB}[n]=\textbf{GL}[n]\times\textbf{S5}$, i.e., $\textbf{GL}[n]\times\textbf{S5}$ is product matching. To prove this, recall 

%\color{blue} I think it is better to incorporate the next proposition into the proof of \cref{Product matching}. \color{black}

We are ready to prove the main result of this section.

\begin{thm} \label{axiomatization_n}
     $\textbf{Grz}[n]\times\textbf{S5}=\textbf{M}^+\textbf{GrzB}[n]$ for each $n \ge 1$.
     %For any $n<\omega$,
\end{thm}
\begin{proof}
    As we already pointed out in \cref{casari_rmk2},  
    $\textbf{M}^+\textbf{GrzB}[n]\subseteq\textbf{Grz}[n]\times\textbf{S5}$. 
    %due to \Cref{casari}. Now 
    For the reverse inclusion, let $\varphi\notin\textbf{M}^+\textbf{GrzB}[n]$. By \crefdefpart{translation}{translation2},  $\varphi^+\notin \textbf{MGLB}[n]$. By \Cref{matching_GLn}, there is a $(\textbf{GL}[n]\times\textbf{S5})$-frame  $\mathfrak{F}\times\mathfrak{G}$ such that $\mathfrak{F}\times\mathfrak{G}\not\models\varphi^+$. Therefore,   $(\mathfrak{F}\times\mathfrak{G})^r\not\models\varphi$ by \Cref{refl_valid}.

    \begin{clm}
        %\color{blue} This can also be numbered (add the claim environment to the preamble) and the proof environment added to its proof (it should read: {\em Proof of the claim}). \color{black} 
        For a $(\textbf{K}\times\textbf{S5})$-frame  $\mathfrak{X}\times\mathfrak{Y}$,  $(\mathfrak{X}\times\mathfrak{Y})^r=\mathfrak{X}^r\times\mathfrak{Y}$. 
        %\color{blue} Why $\frak A$ and $\frak B$? These aren't algebras. \color{black} 
        \label{claim}
    \end{clm}
    \begin{proof*}[Proof of the claim]
       Let $\mathfrak{X}=(X,R_X)$ and $\mathfrak{Y}=(Y,E_Y)$. Then  $\mathfrak{X}\times\mathfrak{Y}=(X\times Y, R, E)$ and $\mathfrak{X}^r\times\mathfrak{Y}=(X\times Y, R',E)$, where $R$, $R'$, and $E$ are defined as in \Cref{products}. It suffices to show that $R' = R^r$. Let $(x,y), (u,v)\in X\times Y$. Then %we have
            \begin{align*}
                (x,y)\mathrel{R'}(u,v) &\iff (x\mathrel{R_X^r} u \text{ and } y=v)\\
                &\iff ((x\mathrel{R_X}u \text{ or } x=u) \text{ and } y=v)\\
                &\iff (x\mathrel{R_X}u \text{ and } y=v) \text{ or } (x=u \text{ and } y=v)\\
                &\iff (x\mathrel{R_X} u \text{ and } y=v) \text{ or } (x,y)=(u,v)\\
                &\iff (x,y)\mathrel{R}(u,v) \text{ or } (x,y)=(u,v)\\
                &\iff (x,y)\mathrel{R^r}(u,v). \qedhere
            \end{align*} 
    \end{proof*}
    
    By \Cref{claim},  $(\mathfrak{F}\times\mathfrak{G})^r=\mathfrak{F}^r\times\mathfrak{G}$. Since $\mathfrak{F}$ is a $\textbf{GL}[n]$-frame, $\mathfrak{F}^r$ is a $\textbf{Grz}[n]$-frame (see the proof of \crefdefpart{frames}{frames1}). Hence, $\varphi$ is refuted in a $(\textbf{Grz}[n]\times\textbf{S5})$-frame. Thus, $\varphi\notin \textbf{Grz}[n]\times\textbf{S5}$, and so $\textbf{Grz}[n]\times\textbf{S5}=\textbf{M}^+\textbf{GrzB}[n]$. 
\end{proof}

%Our final goal will be to prove the analog of the above in the general case, without any finite depth axioms. The first step in doing so is to ensure that $\textbf{M}^+\textbf{GrzB}$ is complete. This we do by proving it has the fmp.

\section{Selective filtration for \textbf{MGrz} and \textbf{MGrzB}} \label{sec4}

In this section, we show that $\textbf{MGrzB}$ admits a form of selective filtration. This will be utilized in \Cref{sec5} to show that \textbf{MGrzB} has the fmp. With appropriate adjustments, 
%to this procedure, 
we also establish the fmp for $\textbf{M}^+\textbf{GrzB}$ and $\textbf{MGLB}$. The selective filtration method employed here is a refinement of the construction developed in \cite{BM}, which itself synthesizes the approaches of \cite{Grefe, mdim, GBm}. %Establishing the fmp for multimodal logics is substantially more intricate than in the unimodal case, and 
The need of the refinement is due to the presence of the monadic Barcan formula, which requires %demands 
that the selected frame satisfies not only the left but also the right commutativity axiom.

%\color{blue}
%Mention that our goal is to prove fmp of MGrzB and MGLB. Then point out that proving fmp of multimodal logics is more involved than that for unimodal logics, give a couple of references. 
%\color{black}

We start by outlining the selection procedure of \cite{BM} for $\textbf{MGrz}$ and then detail the necessary modifications.
While our construction can be performed in any descriptive frame, it is sufficient to only work with canonical models (see, e.g., \cite[Sec.~5.1]{CZ}). %only results were formulated in the more general setting of descriptive frames. However, it is enough to work with canonical models for the logics in question. First, recall the notion of definable subsets.

\begin{defn}
    Let $\mathfrak{F}=(X,R,E)$ be an \textbf{MK}-frame and $v$ a valuation. A subset $U\subseteq X$ is called \emph{definable} in the model $\mathfrak{M} :=(\mathfrak{F},v)$ provided $U=v(\varphi)$ for some $\mathcal{L}_{\Diamond\exists}$-formula $\varphi$.
\end{defn}

%In particular, we work with the definable subsets of the canonical models (see, e.g., \cite[Sec.~5.1]{CZ}) of \textbf{MGrz} and \textbf{MGrzB} to perform our selection. 
We will now restrict our attention to \textbf{MGrz} and \textbf{MGrzB}. We let $\textbf{M}$ be one of these %monadic 
logics, and let $\mathfrak{M}_\textbf{M} :=(\mathfrak{F}_\textbf{M}, v_\textbf{M})$ be the canonical model of \textbf{M}, where $\mathfrak{F}_\textbf{M} :=(X_\textbf{M}, R_\textbf{M},  E_\textbf{M})$ is the canonical frame and $v_\textbf{M}$ the canonical valuation. 
%of \textbf{M}.
%, where $\textbf{M}=\textbf{MGrz} \text{ or } \textbf{MGrzB}$; let . 

%The technique of selecting points originates in \cite{gabbay1970selective, fine_selective}. In particular, the selected points must be \emph{maximal}. For our purposes however, it is not enough to work with maximal points. They must satisfy an additional property, and we deem such points \emph{strongly maximal}. 

Let $U\subseteq X_\textbf{M}$. We recall (see, e.g., \cite[Def.~1.4.9]{esakia2019}) that $x \in U$ is a {\em maximal} point of $U$ provided $x\mathrel{R_\textbf{M}}y$ implies $x=y$ for all $y\in U$. Let ${\textbf{max}}\hspace{2pt}U$ denote the set of maximal points of $U$. It follows from the well-known Fine-Esakia principle (see, e.g., \cite[Cor.~3.5.7]{esakia2019}) that for each definable subset $U\subseteq X_\textbf{M}$ and $x\in U$, we have $R_\textbf{M}[x]\cap\textbf{max}\hspace{2pt}U \neq \varnothing$. In particular, it is not possible to ``exit'' $U$ via $R_\textbf{M}$ from a maximal point and re-enter it; that is, if $x\in\textbf{max}\hspace{2pt}U$, then there are no $y\notin U$ and $z\in U$ such that $x\mathrel{R_\textbf{M}}y\mathrel{R_\textbf{M}}z$ (see \cite[p.~70]{esakia2019}). This will be used in \cref{claim about sim}.

%\color{red}
%\begin{rmk}\label{passive}
%    Recall that in $\mathfrak{F}_\textbf{M}$, for a definable subset $U$ and $x\in\textbf{max}\hspace{2pt}U$, the following situation is forbidden: $x\mathrel{R_\textbf{M}}y\mathrel{R_\textbf{M}}z$, where $y\notin U$ and $z\in U$. In other words, this means that it is not possible to ``exit'' a definable subset via $R_\textbf{M}$ from a maximal point and re-enter it. For details, see \cite[p.~70]{esakia2019}.
%\end{rmk}
%\color{black}

We will work with the stronger notion of maximality, introduced in \cite{BM}.
To define strong maximality, 
%it is convenient to work with the derived relation 
let $Q_\textbf{M}$ be the composition of $R_\textbf{M}$ and $E_\textbf{M}$; that is, %for any $x,y\in X_\textbf{M}$, set 
\begin{equation*}
    x \mathrel{Q_\textbf{M}} y\iff (\exists z\in X_\textbf{M})(x \mathrel{R_\textbf{M}} z \text{ and } z\mathrel{E_\textbf{M}} y).
\end{equation*}

    Clearly $R_\textbf{M}\subseteq Q_\textbf{M}$, but this inclusion is in general strict. 
    %, and we may 
    We will refer to the elements of $Q_\textbf{M} \setminus R_\textbf{M}$ as \emph{proper} $Q$-relations (see \cref{Q-arrow1,Q-arrow2}).
%for some $z\in X_\textbf{M}$. 
%Thus, we have the following definition.

\begin{defn} \cite[Def.~6.13]{BM} Let $U\subseteq X_\textbf{M}$. We call 
%and $x\in U$.
%    \begin{enumerate}        
%        \item \cite[Def.~1.4.9]{esakia2019} 
%        If for all $y\in U$ we have $x\mathrel{R_\textbf{M}}y\implies x=y$, then $x$ is called maximal in $U$. We let ${\textbf{max}}\hspace{2pt}U$ denote the set of points maximal in $U$.
%        
%        \item 
%        If 
$x\in\textbf{max}\hspace{2pt}U$ a {\em strongly maximal} point of $U$ if for any $y\in X_\textbf{M}$,
        \begin{equation*}
            x\mathrel{Q_\textbf{M}}y \text{ and } y\in U \implies  x\mathrel{E_\textbf{M}}y. 
        \end{equation*}
            %then $x$ is called strongly maximal in $U$. 
            We let ${\textbf{smax}}\hspace{2pt}U$ denote the set of  strongly maximal points of $U$.
    %\end{enumerate}
\end{defn}

%\color{blue} Maybe it is worth  explaining the difference and adding an example. \color{black}

The following, taken from \cite[Rmk.~6.14]{BM}, is an illustrative example highlighting the difference between maximal and strongly maximal points. Consider the finite frame $\mathfrak{F}=(X,R,E)$ shown below, where $R$ is represented by the arrows (indicating non-reflexive relations) and circles (indicating reflexive relations), and the $E$-clusters are depicted in blue. Since $X$ is finite, there exists a valuation $v$ under which every subset of $X$ is definable. Let $U=\{a,b,c\}$ (the dotted red curve on the diagram). Clearly $a\in\textbf{max}\hspace{2pt}U$, but $a\notin\textbf{smax}\hspace{2pt}U$ since $a\mathrel{Q}c$, $c\in U$, and $a\mathrel{\cancel{E}}c$. 
%\color{blue} Since $c$ plays no role, it can be removed from the diagram (or delete $d$ and replace it with $c$). \color{black}

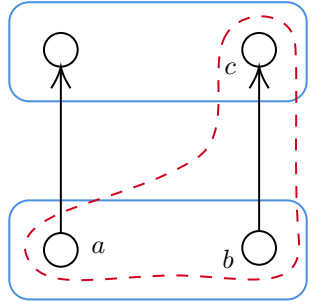
\begin{center}
\tikzset{every picture/.style={line width=0.75pt}} %set default line width to 0.75pt        
\begin{tikzpicture}[x=0.75pt,y=0.75pt,yscale=-1,xscale=1]
%Rounded Rect [id:dp060937344810754634] 
\draw  [color={rgb, 255:red, 74; green, 144; blue, 226 }  ,draw opacity=1 ] (220,81.2) .. controls (220,75.68) and (224.48,71.2) .. (230,71.2) -- (360,71.2) .. controls (365.52,71.2) and (370,75.68) .. (370,81.2) -- (370,111.2) .. controls (370,116.72) and (365.52,121.2) .. (360,121.2) -- (230,121.2) .. controls (224.48,121.2) and (220,116.72) .. (220,111.2) -- cycle ;
%Rounded Rect [id:dp6237168002414025] 
\draw  [color={rgb, 255:red, 74; green, 144; blue, 226 }  ,draw opacity=1 ] (220,181.2) .. controls (220,175.68) and (224.48,171.2) .. (230,171.2) -- (360,171.2) .. controls (365.52,171.2) and (370,175.68) .. (370,181.2) -- (370,211.2) .. controls (370,216.72) and (365.52,221.2) .. (360,221.2) -- (230,221.2) .. controls (224.48,221.2) and (220,216.72) .. (220,211.2) -- cycle ;
%Shape: Circle [id:dp08497447052441975] 
\draw   (237.5,95.2) .. controls (237.5,90.51) and (241.31,86.7) .. (246,86.7) .. controls (250.69,86.7) and (254.5,90.51) .. (254.5,95.2) .. controls (254.5,99.89) and (250.69,103.7) .. (246,103.7) .. controls (241.31,103.7) and (237.5,99.89) .. (237.5,95.2) -- cycle ;
%Shape: Circle [id:dp18205662301626568] 
\draw   (237.5,196.2) .. controls (237.5,191.51) and (241.31,187.7) .. (246,187.7) .. controls (250.69,187.7) and (254.5,191.51) .. (254.5,196.2) .. controls (254.5,200.89) and (250.69,204.7) .. (246,204.7) .. controls (241.31,204.7) and (237.5,200.89) .. (237.5,196.2) -- cycle ;
%Shape: Circle [id:dp1802208523609593] 
\draw   (337.5,95.2) .. controls (337.5,90.51) and (341.31,86.7) .. (346,86.7) .. controls (350.69,86.7) and (354.5,90.51) .. (354.5,95.2) .. controls (354.5,99.89) and (350.69,103.7) .. (346,103.7) .. controls (341.31,103.7) and (337.5,99.89) .. (337.5,95.2) -- cycle ;
%Shape: Circle [id:dp7018613457123046] 
\draw   (337.5,195.2) .. controls (337.5,190.51) and (341.31,186.7) .. (346,186.7) .. controls (350.69,186.7) and (354.5,190.51) .. (354.5,195.2) .. controls (354.5,199.89) and (350.69,203.7) .. (346,203.7) .. controls (341.31,203.7) and (337.5,199.89) .. (337.5,195.2) -- cycle ;
%Shape: Polygon Curved [id:ds5357136822624223] 
\draw  [color={rgb, 255:red, 208; green, 2; blue, 27 }  ,draw opacity=1 ][dash pattern={on 4.5pt off 4.5pt}] (328,90.2) .. controls (335,73.2) and (362,75.2) .. (364,91.2) .. controls (366,107.2) and (363,158.2) .. (366,190.2) .. controls (369,222.2) and (332.09,207.58) .. (297,209.2) .. controls (261.91,210.82) and (239.85,215.59) .. (231,203.2) .. controls (222.15,190.81) and (233,182.2) .. (252,175.2) .. controls (271,168.2) and (295,162.2) .. (314,150.2) .. controls (333,138.2) and (321,107.2) .. (328,90.2) -- cycle ;
%Straight Lines [id:da4997046316950575] 
\draw    (246,187.7) -- (246,105.7) ;
\draw [shift={(246,103.7)}, rotate = 90] [color={rgb, 255:red, 0; green, 0; blue, 0 }  ][line width=0.75]    (10.93,-4.9) .. controls (6.95,-2.3) and (3.31,-0.67) .. (0,0) .. controls (3.31,0.67) and (6.95,2.3) .. (10.93,4.9)   ;
%Straight Lines [id:da8012155573017145] 
\draw    (346,186.7) -- (346,105.7) ;
\draw [shift={(346,103.7)}, rotate = 90] [color={rgb, 255:red, 0; green, 0; blue, 0 }  ][line width=0.75]    (10.93,-4.9) .. controls (6.95,-2.3) and (3.31,-0.67) .. (0,0) .. controls (3.31,0.67) and (6.95,2.3) .. (10.93,4.9)   ;

% Text Node
\draw (260,190.4) node [anchor=north west][inner sep=0.75pt]    {$a$};
% Text Node
\draw (326,194.4) node [anchor=north west][inner sep=0.75pt]    {$b$};
% Text Node
%\draw (256.5,98.6) node [anchor=north west][inner sep=0.75pt]    {$c$};
% Text Node
\draw (327,100.4) node [anchor=north west][inner sep=0.75pt]    {$c$};
\end{tikzpicture}
\end{center}

While the notions of maximality and strong maximality differ in general, they agree for $E$-saturated definable subsets, where we recall that a subset $V$ is {\em $E$-saturated} provided $x \in V$ and $x \mathrel{E} y$ imply $y \in V$. 
%$V=E[V]$, where 
%$$E[V]=\{x\mid y\mathrel{E}x \text{ for some }y\in V\}.$$
Indeed, for such a definable subset $V$, we have $\textbf{smax}\hspace{2pt}V =\textbf{max}\hspace{2pt}V$ (see \cite[Prop.~6.15]{BM}). 
It is worth pointing out that the 
%Note that the 
set $U$ in the example above is not $E$-saturated. 

%\color{red}
%In the canonical frame $\mathfrak{F}_\textbf{M}$, 
We will utilize the following feature of strongly maximal points, which is crucial for our construction. 
%later in the section. satisfy an additional property, which we make note of below.

\begin{thm} Let $U$ be a definable subset of $X_\textbf{M}$. \label{smax_thms}
\begin{enumerate}[label=(\arabic*)]
    \item \cite[Thm.~6.17]{BM} {\label{smax}}
    If $x\in U$, then $Q_\textbf{M}[x]\cap {\textbf{smax}}\hspace{2pt}U \neq \varnothing$.
    \item \cite[Lem.~6.18(2)]{BM} \label{passive_Q}
    If $x\in \textbf{smax}\hspace{2pt} U$, 
    %for some definable subset $U \subseteq X_\textbf{M}$  
$x\mathrel{Q_\textbf{M}}y$, and $y\mathrel{Q_\textbf{M}}x$, then $x\mathrel{E_\textbf{M}}y$.
\end{enumerate}
\end{thm}

%We will utilize the following version of the Fine-Esakia principle.

%\begin{thm}[Fine-Esakia Principle]\cite[Cor.~3.5.7]{esakia2019}{\label{fe_principle}}
%    For each definable subset $U\subseteq X_\textbf{M}$ and $x\in U$, we have $R_\textbf{M}[x]\cap\textbf{max}\hspace{2pt}U \neq \varnothing$.
%\end{thm}

%\begin{thm}\cite[Thm.~6.17]{BM}{\label{smax}}
%    For each definable subset $U\subseteq X_\textbf{M}$ and $x\in U$, we have $Q_\textbf{M}[x]\cap {\textbf{smax}}\hspace{2pt}U \neq \varnothing$. 
%\end{thm}

%\color{blue} Briefly describe the structure of the rest of the section, mentioning that we'll have two subsections and what'll happen in them. \color{black}

%\color{red}
%We will present our observations over two subsections. In the first subsection, we will detail the selective filtration proposed in \cite{BM} for \textbf{MGrz}, following which, we will highlight the necessary modifications in the next subsection.  
%\color{black}

\subsection{Selective filtration for \textbf{MGrz}} \label{subsec: MGrz}

%Suppose $\varphi\notin\textbf{MGrz}$; then $\varphi$ is refuted in 
Let $\mathfrak{M}_\textbf{MGrz} \not\models \varphi$ and let $S$ be the set of subformulae of $\varphi$. We will iteratively construct a new model $\widehat{\mathfrak{M}}$ such that $\widehat{\mathfrak{M}}\not\models\varphi$. 

%In particular, 
By \crefdefpart{smax_thms}{smax}, there is 
%we can find some 
$x\in \textbf{smax}\hspace{2pt}v_\textbf{MGrz}(\neg \varphi)$. To initiate the construction of the model $\widehat{\mathfrak{M}}$, we ``select'' $x$, label it $\widehat{x}$, and proceed by selecting additional points. For $t\in X_\textbf{MGrz}$, we label a selected copy by $\widehat{t}$. 
%Our goal will be to continue this process so that 
%\begin{equation*}    \mathfrak{M}_\textbf{MGrz},t\models\psi\iff \widehat{\mathfrak{M}}, \widehat{t}\models\psi
%\end{equation*}
%for any $\psi\in S$. In particular, to satisfy the above for subformulae of the form $\exists\psi$ or $\Diamond\psi$, we may need to introduce new points in 
New points in the selected model are introduced either as ``witnesses'' of formulae in $S$ of the form $\exists\psi,\Diamond\psi$ or as ``witnesses'' of left commutativity. 

Define
\begin{align*}
    W_t^\exists=\{\exists\psi\in S &\mid \mathfrak{M}_\textbf{MGrz},t\models\exists\psi \text{ but } \mathfrak{M}_\textbf{MGrz},t\not\models\psi\},\\
    W_t^\Diamond=\{\Diamond\psi\in S &\mid \mathfrak{M}_\textbf{MGrz},t\models\Diamond\psi \text{ but } \mathfrak{M}_\textbf{MGrz},t\not\models\psi\}.
\end{align*}

%Therefore, 
The sets $W_t^\exists$ and $W_t^\Diamond$ let us keep track of the subformulae of $\varphi$ that demand the addition of a new witness from a point $\widehat{t}$ in the model $\widehat{\mathfrak{M}}$. A point added to satisfy a subformula of the form $\exists\psi$ will be called an {\em $\exists$-witness}, and a point added to satisfy a subformula of the form $\Diamond\psi$ a {\em $\Diamond$-witness}.

We will drop the subscripts from $\mathfrak{M}_\textbf{MGrz}=(\mathfrak{F}_\textbf{MGrz},v_\textbf{MGrz})$ and $\mathfrak{F}_\textbf{MGrz}=(X_\textbf{MGrz},R_\textbf{MGrz}, E_\textbf{MGrz})$, and simply write $t\models\psi$ to indicate $\mathfrak{M},t\models\psi$. In doing so, no ambiguity will arise since statements about the model $\widehat{\mathfrak{M}}$ will be clearly indicated by the use of hats.
%ted characters.

Define the equivalence relation $\sim_S$ by
\begin{equation*}
    x\sim_S y\iff  x\models \psi \text{ iff } y\models \psi \text{ for any } \psi\in S. 
\end{equation*}

\begin{enumerate}[label={(\arabic*)}]
    \item \textbf{Initial step}
    
        Select $x\in \textbf{smax}\hspace{2pt}v(\neg \varphi)$ and let 
        \begin{equation*}
            X_0=\{\widehat{x}\}, R_0=X_0^2, E_0=X_0^2.
        \end{equation*}

        Set $\mathfrak{F}_0=(X_0,R_0,E_0)$ and note that $\mathfrak{F}_0$ is a partially ordered \textbf{MS4}-frame.

        Suppose a partially ordered \textbf{MS4}-frame $\mathfrak{F}_{k-1}=(X_{k-1},R_{k-1},E_{k-1})$ has already been constructed. We describe the construction of a partially ordered \textbf{MS4}-frame $\mathfrak{F}_k=(X_k,R_k,E_k)$.
        
    \item $\exists$ \textbf{step}: Selecting $\exists$-witnesses.

        Set $X_k^\exists=X_{k-1}, R_k^\exists=R_{k-1}, E_k^\exists=E_{k-1}$. 

        Suppose $\widehat{t}\in X_{k-1}$ and $\exists\psi\in W_t^\exists$. If there is already some $\widehat{u}\in X_k^\exists$ such that $\widehat{t}\mathrel{E_k^\exists}\widehat{u}$ and $u\models\psi$, there is nothing to do. Otherwise, we choose a witness using the following lemma.

        \begin{lem}\cite[Lem.~7.2]{BM}{\label{horz}}
            Let $z\in{\textbf{smax}}\hspace{2pt}U$ for some definable subset $U$. If $\exists\psi\in W_z^{\exists}$, then there is $y \in {\textbf{smax}}\hspace{2pt}(E[U]\cap v(\psi))$ such that $y \mathrel{E} z$. 
        \end{lem}

        Since $\widehat{t}\in X_{k-1}$, and we only select strongly maximal points at each stage of the construction, it follows that $t\in\textbf{smax}\hspace{2pt} U$ for some definable subset $U$. Therefore, by \Cref{horz}, there is some 
        $w\in\textbf{smax}\hspace{2pt}(E[U]\cap v(\psi))$ such that $w\mathrel{E}t$ (in particular, note that $w\models\psi$ and that $E[U]\cap v(\psi)$ is definable).

        Add $(\widehat{t},\widehat{w})$ to $E_k^\exists$ and generate the least equivalence relation. Add $(\widehat{w},\widehat{w})$ to $R_k^\exists$. Repeat this process until all points have witnesses for subformulae of the form $\exists\psi$.

    \item $\Diamond$ \textbf{step}: Selecting $\Diamond$-witnesses.

    Set $X_k^\Diamond=X_k^\exists, R_k^\Diamond=R_k^\exists, E_k^\Diamond=E_k^\exists$.

    Suppose $\widehat{y}\in X_k^\exists$ and $\Diamond\psi\in W_y^\Diamond$. If there is already some $\widehat{z}\in X_k^\Diamond$ such that $\widehat{y}\mathrel{R_k^\Diamond} \widehat{z}$ and $z\models \psi$, there is nothing to do. Otherwise, we use the following lemma to find a witness.

    \begin{lem}\cite[Lem.~7.3]{BM}{\label{vert}}
    For $\Diamond\psi\in W_y^{\Diamond}$ let 
    \begin{equation*}
            A = v(\Diamond\psi)\cap\bigcap\bigl\{ v(\neg\Diamond\alpha)\mid \Diamond\alpha\in S,\,  y\not\models \Diamond\alpha \bigr\}.
        \end{equation*}
    Then there is $z\in X$ such that the following conditions are satisfied:
    \begin{enumerate}[label={(\arabic*)}]
        \item $y\mathrel{Q}z$, $z\neq y$, and $z\in\textbf{smax}\hspace{2pt}A\cap\textbf{max}\hspace{2pt}v(\psi)$; \label{vert3}
        
        \item $y\not\models\Diamond\alpha$ implies $z\not\models\alpha$ for each $\Diamond\alpha\in S$; \label{vert2}
        
        \item If $y\mathrel{E}z$, then there is $u\in {\textbf{smax}}\hspace{2pt}A \cap {\textbf{max}}\hspace{2pt}v(\psi)$ such that $y\mathrel{R}u$ and $u\mathrel{E}y$; \label{vert4}

        \item If $z\mathrel{Q}t$ and $z \mathrel{\cancel{E}} t$, then $t\mathrel{\cancel{\sim}_S}y$. \label{vert6}
    \end{enumerate}
    \end{lem}

    %\color{blue} The above lemma is technical and hard to follow. Do we use all the items in what follows? If not, only leave the ones that are essential. That would make it easier to comprehend. \color{black} 

    \textbf{Horizontal step}: Selecting $\Diamond$-witnesses in $E$-clusters.

    We will select points $u$ satisfying
    \begin{equation}
        y\mathrel{R}u, y\mathrel{E}u, \text{ and }  u\in\textbf{smax}\hspace{2pt} A\cap\textbf{max}\hspace{2pt} v(\psi).  \label{hw}\tag{H}
    \end{equation}

    \begin{enumerate} 
        \item[(h1)] Suppose there is already some $\widehat{u}\in X_k^\Diamond$ such that \labelcref{hw} holds for $u$. 

        \item[(h2)] Then add $(\widehat{y},\widehat{u})$ to $R_k^\Diamond$ and take the reflexive-transitive closure. Further, add $(\widehat{y},\widehat{u})$ to $E_k^\Diamond$ and generate the least equivalence relation. \label{h2}

        \item[(h3)] If there is no such $\widehat{u}$, then by \crefdefpart{vert}{vert3} 
        %\color{blue} (the labeling needs fixing; there should be no 3 just a,c (or simply c as a is part of c); same problem where the lemma is referenced below) \color{black} 
        there is some $z\in \textbf{smax}\hspace{2pt}A\cap\textbf{max}\hspace{2pt}v(\psi)$ such that $y\mathrel{Q}z$. If additionally we have $y\mathrel{E}z$, then by \crefdefpart{vert}{vert4} there is 
        %a point 
        $u$ satisfying \labelcref{hw}. We add $\widehat{u}$ to $X_k^\Diamond$ and introduce relations as in (\hyperref[h2]{h2}).
    \end{enumerate}

    Repeat the above steps until all such $\Diamond$-witnesses for $\widehat{y}$ from $E[y]$ are added. Proceed to the next step.

    \begin{rmk}\label{horz_rmk}
        %We point out that 
        This is the only step in the construction where a ``horizontal'' $R_k^\Diamond$-arrow is drawn, i.e., relations of the form $\widehat{y}\mathrel{R_k^\Diamond}\widehat{u}$ and $\widehat{y}\mathrel{E_k^\Diamond}\widehat{u}$ are introduced. 
    \end{rmk}
    
    \textbf{Vertical step}: Selecting $\Diamond$-witnesses outside $E$-clusters.

     We will select points $z$ 
     %\color{blue} (why not call these $u$ as in the horizontal step?) \color{black} 
     satisfying
    \begin{equation}
        y\mathrel{Q}z, y\mathrel{\cancel{E}}z, \text{ and }  z\in\textbf{smax}\hspace{2pt} A\cap\textbf{max}\hspace{2pt} v(\psi).  \label{vw}\tag{V}
    \end{equation}

    \begin{enumerate} 
        \item[(v1)] Suppose there is already some $\widehat{z}\in X_k^\Diamond$ such that \labelcref{vw} holds for $z$. Then add $(\widehat{y},\widehat{z})$ to $R_k^\Diamond$ and take the reflexive-transitive closure. 

        \item[(v2)] If there is no such $\widehat{z}$, then by \crefdefpart{vert}{vert3} there is some $z\in \textbf{smax}\hspace{2pt}A\cap\textbf{max}\hspace{2pt}v(\psi)$ such that $y\mathrel{Q}z$. Without loss of generality we may assume that
        %Then we will also have 
    $y\mathrel{\cancel{E}}z$ since the case $y\mathrel{E}z$ is already considered in the horizontal step.  Clearly $z$ satisfies \labelcref{vw}. We add $\widehat{z}$ to $X_k^\Diamond$, add $(\widehat{y},\widehat{z})$ to $R_k^\Diamond$ and take the reflexive-transitive closure, and add $(\widehat{z},\widehat{z})$ to $E_k^\Diamond$. 

        \item[(v3)] If $z\mathrel{E}t$ and $\widehat{t}$ has been previously introduced, then add $(\widehat{z},\widehat{t})$ to $E_k^\Diamond$ and generate the least equivalence relation.
    \end{enumerate}

    Repeat the above steps until all $\Diamond$-witnesses for $\widehat{y}$ occurring outside $E[y]$ are added. 

    Repeat the horizontal and vertical steps for all points in $X_k^\exists$.

    \begin{rmk}\label{Q-arrow1}
        %Note that 
        The vertical step may involve turning a proper $Q$-arrow into an $R_k^\Diamond$-arrow. In our final construction for \textbf{MGrzB}, this will be one of the two instances where this occurs (for the other instance, see \Cref{Q-arrow2}).
    \end{rmk}

%\color{blue} Combine the next two lemmas into one. \color{black}

    \begin{lem} \label{EQ-lemma}
        \cite[Lem.~7.4]{BM}
        For $\widehat{u},\widehat{w}\in X_k^\Diamond$,
        \begin{enumerate}
            \item  $\widehat{u}\mathrel{E_k^\Diamond}\widehat{w}$ iff $u\mathrel{E}w$. {\label{E-lemma}}

            \item $\widehat{u}\mathrel{R_k^\Diamond}\widehat{w}$ implies $u\mathrel{Q}w$. \label{Q-lemma}
        \end{enumerate}
    \end{lem}

    \item \textbf{Commutativity step}

    In this final step, we close the selected points under left commutativity. The added points also require closure under left commutativity, 
    %which results in the addition of new points, which again demand witnesses for left commutativity, 
    warranting further applications of this step. However, the process terminates in finitely many steps. 
    %We stratify this process as follows.

    Set $X_k^0=X_k^\Diamond, R_k^0=R_k^\Diamond, E_k^0=E_k^\Diamond$. Suppose $X_k^{j-1}, R_k^{j-1}, E_k^{j-1}$ have already been defined. Let $X_k^j=X_k^{j-1},R_k^j=R_k^{j-1},E_k^j=E_k^{j-1}$. We update these sets to ensure that left commutativity is satisfied for points in $X_{j-1}$, i.e., the required witnesses will be included in $X_k^j$ (along with the necessary $R_k^j$ and $E_k^j$ relations). At each stage indexed by $j$, we will ensure that an analog of \Cref{EQ-lemma} is applicable (see \Cref{lemma_lc}). 
    
    Suppose $\widehat{t}\mathrel{E_k^{j-1}}\widehat{u}$ and $\widehat{u}\mathrel{R_k^{j-1}}\widehat{w}$. We need to find some $\widehat{s}$ such that $\widehat{t}\mathrel{R_k^j}\widehat{s}$ and $\widehat{s}\mathrel{E_k^j}\widehat{w}$. In case $\widehat{u}\mathrel{E_k^{j-1}}\widehat{w}$, the choice of $\widehat{s}$ is trivial since we can simply take $\widehat{s}=\widehat{t}$. Hence, suppose $\widehat{u}\mathrel{\cancel{E_k^{j-1}}}\widehat{w}$. Then we can find 
    %a witness 
    $\widehat{s}$ using the following: 
    %lemma.

    \begin{lem}\cite[Lem.~7.5]{BM}\label{l-comm}
        Let $t\mathrel{Q}w$ and $w\in\textbf{smax}\hspace{2pt}U$ for some definable set $U$. Then there exists $s\in\textbf{smax}\hspace{2pt}E[U]$ such that $t\mathrel{R}s$ and $s\mathrel{E}w$. %\color{blue} We should mention earlier that $\textbf{smax}U=\textbf{max}U$ for $E$-saturated $U$. \color{black}
    \end{lem}

    Since $\widehat{t}\mathrel{E_k^{j-1}}\widehat{u}$ and $\widehat{u}\mathrel{R_k^{j-1}}\widehat{w}$, \Cref{lemma_lc} (or \Cref{EQ-lemma} for $j=1$) gives that $t\mathrel{E}u$ and $u\mathrel{Q}w$. Thus, $t\mathrel{Q}w$, where $w\in\textbf{smax}\hspace{2pt}U$ for some definable subset $U$ (since every selected point is strongly maximal in some definable subset). Therefore, \Cref{l-comm} applies, by which there is $s\in\textbf{smax}\hspace{2pt}E[U]$ such that $t\mathrel{R}s$ and $s\mathrel{E}w$.

    If $\widehat{s}$ is not already present, add it to $X_k^j$. Also add $(\widehat{t},\widehat{s})$ to $R_k^j$ and take the reflexive-transitive closure. In addition, add $(\widehat{s},\widehat{w})$ to $E_k^j$ and generate the least equivalence relation. Repeat this process for every instance of $\widehat{t}\mathrel{E_k^{j-1}}\widehat{u}$ and $\widehat{u}\mathrel{R_k^{j-1}}\widehat{w}$ until left commutativity is satisfied for all the points in $X_k^{j-1}$. 

    \begin{lem}\cite[Lem.~7.6]{BM}{\label{lemma_lc}}
    For $\widehat{u},\widehat{w}\in X_k^j$,
    \begin{enumerate}[label={(\arabic*)}]
        \item $\widehat{u}\mathrel{E_k^j}\widehat{w}$ iff $u\mathrel{E}w$. \label{E_lc}

        \item $\widehat{u}\mathrel{R_k^j}\widehat{w}$ implies $u\mathrel{Q}w$. \label{Q_lc}
    \end{enumerate}
    \end{lem}
    
    Set $X_k=\bigcup_{j<\omega}X_k^j, R_k=\bigcup_{j<\omega}R_k^j, E_k=\bigcup_{j<\omega}E_k^j$. We have thus constructed $\mathfrak{F}_k=(X_k,R_k,E_k)$.

    Note that there is $j$ such that $X_k^j=X_k^{j+l}$ for all $l<\omega$, i.e., there is a stage after which no more points are added, and so $X_k$ is finite (see \cite[Lem.~7.8]{BM}). In contrast, as we will see in \Cref{subsec2}, the set $X_k$ need not be finite when right commutativity is also included in the construction. 
\color{black}
\end{enumerate}

We now define $\widehat{\mathfrak{F}}=(\widehat{X},\widehat{R},\widehat{E})$ by setting %as follows
\begin{equation*}
    \widehat{X}=\bigcup_{k<\omega}X_k, \quad \widehat{R}=\bigcup_{k<\omega}R_k, \quad \widehat{E}=\bigcup_{k<\omega}E_k.
\end{equation*}

%\color{blue} What follows should probably be in the theorem environment. \color{black}
%In \cite[Sec.~7]{BM} we proved the following facts about $\widehat{\mathfrak{F}}$:

\begin{thm} 
\begin{enumerate}[label={(\arabic*)}]
    \item[] 
    \item \cite[Thm.~7.19]{BM} $\widehat{\mathfrak{F}}$ is a finite partially ordered \textbf{MS4}-frame. \label{partialorder}

    \item \cite[Lem.~7.10]{BM} There is a valuation $\widehat{v}$ on $\widehat{\mathfrak{F}}$ such that
    %model $\widehat{\mathfrak{M}}=(\widehat{\mathfrak{F}},\widehat{v})$ such that 
    for each $\widehat{t}\in \widehat{X}$ and $\psi\in S$ we have $$\widehat{t}\models\psi \iff t\models\psi.$$ \label{model}

    %\item \color{blue} Combine (1) and (3). \color{black} \label{finiteframe}
\end{enumerate}
\end{thm}

In the following subsection, we adjust the above construction for \textbf{MGrzB}. 
%Between \cite{BM} and the construction proposed in the following subsection, 
The key difference is that in addition to left commutativity,
%in the commuativity step. Indeed, we only required left commmutativity in \cite{BM}. However, forcing the construction to 
we also need to satisfy right commutativity. Due to this, the procedure may not terminate in finitely many steps (closing under left commutativity necessitates closing under right commutatvity and vice-versa; thus, these two processes alternate and trigger each other). Therefore, the selected frame $\widehat{\mathfrak{F}}$ 
%obtained through this construction 
will be countable, but not necessarily finite. 
%in general. 
However, 
%as we will see, 
the frame will 
%refute the formula $\varphi$ (i.e., it satifies \labelcref{model} above) and will 
be an \textbf{MGrzB}-frame with a uniform bound %$n$ 
on its depth. 

%\color{blue} Explain that this is needed to also satisfy the right commutativity, and that this will be achieved by sacrificing the finiteness of $\widehat{\mathfrak F}$. \color{black} 

\subsection{Selective filtration for \textbf{MGrzB}} \label{subsec2}

Let $\mathfrak{M}=(\mathfrak{F},v)$ be the canonical model of \textbf{MGrzB}. Since the monadic Barcan formula and its converse are Sahlqvist, they are canonical (see, e.g., \cite[p.~223]{mdim}). Thus, the frame $\mathfrak{F}$ satisfies both the left and right commutativity. 

We perform the same construction as in \cref{subsec: MGrz}, but  modify the commutativity step accordingly. 

\textbf{Commutativity step}

Set $X_k^0=X_k^\Diamond, R_k^0=R_k^\Diamond, E_k^0=E_k^\Diamond$. For each $j<\omega$, we define 
%the sets 
$X_k^j, R_k^j, E_k^j$ as follows: if $j-1$ is even, apply the step \labelcref{LC}; if $j-1$ is odd, apply the step \labelcref{RC}. %In other words, 
After every even indexed step, the constructed frame will satisfy left commutativity for the points added previously, and after every odd indexed step,  
right commutativity. At each stage, we will ensure that \Cref{EQ-lemma} is still satisfied, i.e., that $\widehat{u}\mathrel{E_k^j}\widehat{w}$ iff $u\mathrel{E}w$, and $\widehat{u}\mathrel{R_k^j}\widehat{w}$ implies $u\mathrel{Q}w$.

Because the construction 
%now 
involves successive applications of 
%the steps 
\labelcref{LC} and \labelcref{RC}, each application has its own analog of \Cref{EQ-lemma}. At each odd $j$, we have \Cref{EQ-lemma_lc}; at each even $1<j$, we have \Cref{EQ-lemma_rc};
%similarly 
%for every even $1<j$. 
%Note that  
and at $j=0$, \Cref{EQ-lemma} applies. 

\begin{enumerate}[label={(LC)}]
    \item \textbf{Left commutativity} \label{LC}

    Set $X_k^j=X_k^{j-1}, R_k^j=R_k^{j-1}, E_k^j=E_k^{j-1}$.

    Suppose $\widehat{t}\mathrel{E_k^{j-1}}\widehat{u}$ and $\widehat{u}\mathrel{R_k^{j-1}}\widehat{w}$. By \Cref{l-comm,EQ-lemma_rc} (for $j-1$), 
    %\color{blue} (the latter is not proved yet), \color{black} 
    we may find 
    a witness $\widehat{s}$ for left commutativity, as described in the previous subsection. We introduce $\widehat{s}$ to $X_k^j$, add $(\widehat{t},\widehat{s})$ to $R_k^j$ and take the reflexive-transitive closure. We also add $(\widehat{s},\widehat{w})$ to $E_k^j$ and generate the least equivalence relation. This is then repeated until left commutativity is satisfied for all the points in $X_k^{j-1}$. 

    The proof of the following lemma is similar to that of \cite[Lem.~7.6]{BM}.
    
    \begin{lem}{\label{EQ-lemma_lc}}
    For $\widehat{u},\widehat{w}\in X_k^j$,
    \begin{enumerate}[label={(\arabic*)}]
        \item $\widehat{u}\mathrel{E_k^j}\widehat{w}$ iff $u\mathrel{E}w$. \label{E-lemma_lc}

        \item $\widehat{u}\mathrel{R_k^j}\widehat{w}$ implies $u\mathrel{Q}w$. \label{Q-lemma_lc}
    \end{enumerate}
    \end{lem}
    \begin{proof}
        \labelcref{E-lemma_lc} The implication $\widehat{u}\mathrel{E_k^j}\widehat{w}\implies u\mathrel{E}w$ is clear since we introduce $E_k^j$-relations between points only if their original copies are $E$-related. Conversely, suppose $u\mathrel{E}w$. Since $\widehat{u},\widehat{w}\in X_k^j$, 
        there are $\widehat{t}, \widehat{s}\in X_k^{j-1}$ such that $\widehat{w}\mathrel{E_k^j}\widehat{t}$ and $\widehat{u}\mathrel{E_k^j}\widehat{s}$ (because each point that may be introduced via \labelcref{LC} must involve creating an $E_k^j$-relation for left commutativity). Thus, by the forward implication, 
        $t\mathrel{E}w$ and $u\mathrel{E}s$, and hence $t\mathrel{E}s$. By \Cref{EQ-lemma_rc} (or \Cref{EQ-lemma}), we  deduce that $\widehat{s}\mathrel{E_k^{j-1}}\widehat{t}$, and so $\widehat{t}\mathrel{E_k^j}\widehat{s}$. %Hence, we have 
        Thus, $\widehat{u}\mathrel{E_k^j}\widehat{w}$.

        \labelcref{Q-lemma_lc} 
        This is immediate since a $R_k^j$-relation is introduced between points only if their original copies are $Q$-related.
        %, and the result follows. 
    \end{proof}
\end{enumerate}    

%\color{blue} Need to say something about the proof of the above lemma. Is it the same as in the previous subsection? \color{black}

\begin{enumerate}[label={(RC)}]
    \item \textbf{Right commutativity} \label{RC}

    Set $X_k^j=X_k^{j-1}, R_k^j=R_k^{j-1}, E_k^j=E_k^{j-1}$.
    
    Suppose $\widehat{t}\mathrel{R_k^{j-1}}\widehat{w}$ and $\widehat{w}\mathrel{E_k^{j-1}}\widehat{u}$. We need to find $\widehat{s}$ such that $\widehat{t}\mathrel{E_k^j}\widehat{s}$ and $\widehat{s}\mathrel{R_k^j}\widehat{u}$. If $\widehat{t}\mathrel{E_k^{j-1}}\widehat{w}$, there is nothing to do since we can take $\widehat{s}=\widehat{u}$. Hence, suppose that $\widehat{t}\mathrel{\cancel{E_k^{j-1}}}\widehat{w}$. We will find $\widehat{s}$ using the following:

    \begin{lem}\label{r-comm}
        Let $t\in\textbf{smax}\hspace{2pt}U$ for some definable set $U$. If $t\mathrel{Q}u$, then there is $s\in\textbf{smax}\hspace{2pt}(E[U]\cap V)$ such that $s\mathrel{Q}u$ and $s\mathrel{E}t$, where
        \begin{equation*}
            V=\bigcap\{v(\Diamond\psi)\mid \Diamond\psi\in S, u\models\Diamond\psi\}.
        \end{equation*}
    \end{lem}
    \begin{proof}
    Since $t\mathrel{Q}u$, there is $a$ with $t\mathrel{R}a$ and $a\mathrel{E}u$. Because $\mathfrak{F}$ satisfies right commutativity, there is $b$ with $t\mathrel{E}b$ and $b\mathrel{R}u$. The above intersection defining $V$ is finite, and thus $V$ is a definable set. Moreover, $b\in V \cap E[U]$ 
    %and from $t\in U$ and $t\mathrel{E} b$ it follows that $b\in E[U]$. Thus, $b\in E[U]\cap V$. 
    and $E[U]\cap V$ is a definable set. By \crefdefpart{smax_thms}{smax}, we can choose $s\in Q[b]\cap \textbf{smax}\hspace{2pt}(E[U]\cap V)$. We claim that $s$ is the desired point. We first show that $s\mathrel{E}t$. Because $t\mathrel{E}b$ and $b\mathrel{Q}s$, we have $t\mathrel{Q}s$. Since $s\in E[U]$, there is $c \in U$ with $c\mathrel{E}s$. Therefore, $t\mathrel{Q}c$ and $c\in U$. Because $t\in\textbf{smax}\hspace{2pt}U$, we must have $c\mathrel{E}t$, which yields $s\mathrel{E}t$. Now observe that $t\mathrel{Q}u$ and $s\mathrel{E}t$. Thus, 
    %. Thus, we obtain 
    $s\mathrel{Q}u$, concluding the proof.
    %\color{black}
    %For this it is sufficient to show that $s\mathrel{E}t$ 
    %and $s\mathrel{Q}u$. It is further enough to show that $s\mathrel{E}t$, 
    %since $s\mathrel{Q}u$ follows from 
    %the fact that 
    %$t\mathrel{Q}u$ \color{blue} (why?). \color{red} we can now delete this sentence. \color{black} 
    \end{proof}

    Since $\widehat{t}\mathrel{R_k^{j-1}}\widehat{w}$ and $\widehat{w}\mathrel{E_k^{j-1}}\widehat{u}$, we have $t\mathrel{Q}w$ and $w\mathrel{E}u$ by \Cref{EQ-lemma_lc}. This yields $t\mathrel{Q}u$, and by construction, $t\in\textbf{smax}\hspace{2pt}U$ for some definable subset $U$. By \Cref{r-comm}, there is $s\in\textbf{smax}\hspace{2pt}(E[U]\cap V)$ with $s\mathrel{Q}u$ and $s\mathrel{E}t$, where $V$ is defined as above.

    If $\widehat{s}$ is not already there,
    add $\widehat{s}$ to $X_k^j$. Also, add $(\widehat{s},\widehat{u})$ to $R_k^j$ and take the reflexive-transitive closure, and add $(\widehat{t},\widehat{s})$ to $E_k^j$ and generate the least equivalence relation. Repeat until right commutativity is satisfied for all the points in $X_k^{j-1}$.

    \begin{rmk}\label{Q-arrow2}
        %Once again, note that 
        In this step, a proper $Q$-arrow is turned into an $R_k^j$-arrow. This is the second instance in the construction where this occurs (see \Cref{Q-arrow1} for the first instance).
    \end{rmk}

\begin{lem}{\label{EQ-lemma_rc}}
    For $\widehat{u},\widehat{w}\in X_k^j$,
    \begin{enumerate}
        \item $\widehat{u}\mathrel{E_k^j}\widehat{w}$ iff $u\mathrel{E}w$. \label{E-lemma_rc}

        \item $\widehat{u}\mathrel{R_k^j}\widehat{w}$ implies $u\mathrel{Q}w$. \label{Q-lemma_rc}
    \end{enumerate}
    \end{lem}
    \begin{proof}
        The proof is analogous to that of \Cref{EQ-lemma_lc} (but uses 
        %. The claim follows from \Cref{EQ-lemma_lc} for $j-1$. \color{red} Done \color{blue} (maybe point out that the proof uses 
        \Cref{EQ-lemma_lc} instead of \Cref{EQ-lemma_rc}).
    \end{proof}

\end{enumerate}

%\color{blue} How many steps are there? This part needs more detail. \color{black}

Set 
\begin{equation*}
    X_k=\bigcup_{j<\omega}X_k^j, \quad R_k=\bigcup_{j<\omega}R_k^j, \quad
    E_k=\bigcup_{j<\omega}E_k^j,
\end{equation*}
%We have thus constructed 
and let $\mathfrak{F}_k=(X_k,R_k,E_k)$.

%\color{blue}
%Add diagram. This might be a good spot.
%\color{black}
%\color{red}
%Will add a detailed description soon.
%\color{blue}
%Let's not wait on this.
%\color{black}

For the reader's convenience, we depict the various stages of the construction in the diagram below. As usual, circles indicate reflexive points, arrows represent non-reflexive relations in $R_k$, and $E_k$-clusters are in blue. At each stage, the newly added arrows and points are shown in red, while the arrows and points added in the preceding steps are gray and black, respectively.

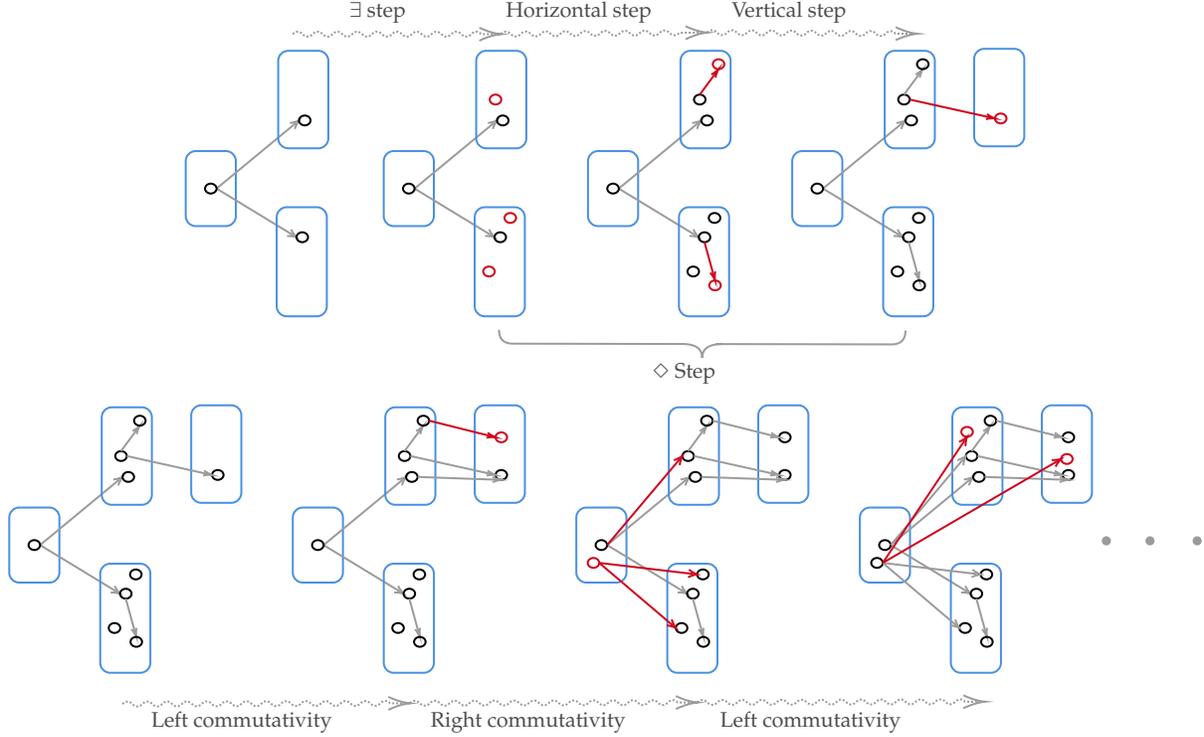
\begin{center}

\tikzset{every picture/.style={line width=0.75pt}} %set default line width to 0.75pt        


\end{center}

Note the alternating nature of the left and right commutativity steps. Closing under left commutativity necessitates closing under right commutativity, and vice-versa. Therefore, while our starting point is a single point in the initial (left-most) cluster, after closing under both left and right commutativity, the initial cluster may acquire new points. These updates to previously ``stabilized'' $E_k$-clusters may result in a potentially infinite frame.

We next outline some
%now make note of the salient 
features of the construction. This 
%We will 
requires several results from \cite[Sec.~7]{BM}. If a proof 
is unaffected by the modification, we simply state the result. 
%For cases where 
If the above adjustment does alter the proof, we supply all the necessary details. The proof of the next lemma is similar to that of \cite[Lem.~7.7]{BM}. 

\begin{lem}\label{rmk}
Suppose $\widehat{u},\widehat{w}\in X_k$.
\begin{enumerate}[label={(\arabic*)}]
    \item $\widehat{u}\mathrel{E_k}\widehat{w}$ iff $u\mathrel{E}w$. \label{rmk1}

    \item $\widehat{u}\mathrel{R_k}\widehat{w}$ implies $u\mathrel{Q}w$. \label{rmk2}
    
    \item $\widehat{u}\mathrel{R_k}\widehat{w}$ implies ($u\mathrel{R}w$ and $u\mathrel{E}w$) or ($u\mathrel{Q}w$ and $u \mathrel{\cancel{E}} w$). \label{rmk3}
    
    %\item $\widehat{u}\mathrel{Q_k}\widehat{w}$ implies $u\mathrel{Q}w$. \label{rmk4}
    
    \item $u=w$ iff $\widehat{u}=\widehat{w}$. \label{rmk4}
    %\color{blue} It is not explained what this means. \color{black} 
\end{enumerate}
\end{lem}

\begin{proof}
    It is enough to show \labelcref{rmk3} since \labelcref{rmk1} and \labelcref{rmk2} are immediate from \Cref{EQ-lemma_lc,EQ-lemma_rc}, and \labelcref{rmk4} is clear (we introduce a hatted point to $X_k$ only if it has not been added previously, and distinct points have distinct hatted versions in $X_k$). 
    %Hence, we show \labelcref{rmk3}. 
    Let $\widehat{u}\mathrel{R_k}\widehat{w}$. %Suppose additionally that 
    If $\widehat{u}\mathrel{\cancel{E_k}}\widehat{w}$, applying \labelcref{rmk1} and \labelcref{rmk2} yields $u\mathrel{Q}w$ and $u\mathrel{\cancel{E}}w$. Suppose that $\widehat{u}\mathrel{E_k}\widehat{w}$. Then $u\mathrel{E}w$ by \labelcref{rmk1}. Because we take the reflexive-transitive closure at each step, we have $\widehat{u}\mathrel{R_k}\widehat{x}_1\mathrel{R_k}\dots\widehat{x}_{n-1}\mathrel{R_k}\widehat{w}$, where each $R_k$-arrow is introduced in one of the steps of the construction. By \labelcref{rmk2}, $u\mathrel{Q}x_1\mathrel{Q}\dots x_{n-1}\mathrel{Q}w$. We claim that $u\mathrel{E}x_i$ for each $1\leq i\leq n-1$. We have $u\mathrel{Q}x_i\mathrel{Q}w\mathrel{Q}u$, where $w\mathrel{Q}u$ follows from $w\mathrel{E}u$. Therefore, 
    %we have 
    $u\mathrel{Q}x_i\mathrel{Q}u$, and since $u$ is strongly maximal, $u\mathrel{E}x_i$ %($1\leq i\leq n-1$) 
    (see \crefdefpart{smax_thms}{passive_Q}). 
    %Together, we have 
    Thus, $\widehat{u}\mathrel{R_k}\widehat{x}_1\mathrel{R_k}\dots\widehat{x}_{n-1}\mathrel{R_k}\widehat{w}$ and  $\widehat{u}\mathrel{E_k}\widehat{x}_1\mathrel{E_k}\dots\widehat{x}_{n-1}\mathrel{E_k}\widehat{w}$, where the latter follows from \labelcref{rmk1}. Applying \labelcref{hw} and \Cref{horz_rmk}, 
    %and \labelcref{hw}, 
    we obtain 
    $u\mathrel{R}x_1\mathrel{R}\dots x_{n-1}\mathrel{R}w$, and hence $u\mathrel{R}w$.
\end{proof}

%\color{blue} Point out that the proof is similar to \cite[Lem.~7.5]{BM}, but add it since the construction is technically different. \color{black}

%\color{blue}
%In \cref{vert} it was a,b,c.. In the above lemma it is 1,2,3.. Make it consistent. Also, are all items necessary? Comment out the ones that aren't. 
%\color{black}

%\color{blue} Some of this already needs to be mentioned up above. \color{black}
%\color{red}
%I have moved the content here to the beginning of the subsection.
%\color{black}
%\color{blue} It is probably enough to mention that the constructed frame will have finite $R$-depth. Local tabularity will be employed in the next section and can be mentioned there. \color{black}

%\color{blue} Because B is present, I'd add a short proof to the next lemma. That it is not necessarily finite is already mentioned above. If need be, we can add a remark after the proof pointing this out again. \color{black}
The following is a weakened version of \cite[Lem.~7.8]{BM}. In particular, the frame constructed in \cite[Lem.~7.8]{BM} is finite, while for the modified construction we may only claim that it is countable. 
%$\mathfrak{F}_k$ is finite , but that may not be the case for the modified construction. However, the following result is still true, and the proof follows that of \cite[Lem.~7.6]{BM}. 

\begin{lem}{\label{partial}}
    $\mathfrak{F}_k$ is a countable partially ordered \textbf{MS4B}-frame for each $k<\omega$.
\end{lem}
\begin{proof}
    That $\mathfrak{F}_k$ is partially ordered with respect to $R_k$ is a consequence of the fact that we only choose strongly maximal points. The details can be found in \cite[Lem.~7.8]{BM}. Additionally, each $X_k^j$ is countable %\color{blue} (even finite) \color{black} 
    for each $j<\omega$, and so $X_k$ is countable as well.
    
    We claim that $\mathfrak{F}_k$ satisfies both left and right commutativity, and hence is an $\textbf{MS4B}$-frame. For 
    %We will prove that 
    left commutativity, 
    %is satisfied; the proof for right commutativity is entirely analogous. Therefore, 
    let $\widehat{t}\mathrel{E_k}\widehat{u}$ and $\widehat{u}\mathrel{R_k}\widehat{w}$. By the definition of $\mathfrak{F}_k$, there exists $j<\omega$ such that $\widehat{t}\mathrel{E_k^j}\widehat{u}$ and $\widehat{u}\mathrel{R_k^j}\widehat{w}$. 
    %Then 
    By construction, there exists $\widehat{s}\in X_k^{j+2}$ such that $\widehat{t}\mathrel{R_k^{j+2}}\widehat{s}$ and $\widehat{s}\mathrel{E_k^{j+2}}\widehat{w}$. Hence, $\widehat{t}\mathrel{R_k}\widehat{s}$ and $\widehat{s}\mathrel{E_k}\widehat{w}$.
    The proof for right commutativity is 
    %entirely 
    analogous. 
\end{proof}

%As in \cref{subsec: MGrz},
%before, 
Now, as in \cref{subsec: MGrz}, define 
%a frame 
$\widehat{\mathfrak{F}}=(\widehat{X},\widehat{R},\widehat{E})$ by setting %taking
\begin{equation*}
    \widehat{X}=\bigcup_{k<\omega}X_k, \quad \widehat{R}=\bigcup_{k<\omega}R_k, \quad \widehat{E}=\bigcup_{k<\omega}E_k.
\end{equation*}

As a consequence of \Cref{partial}, 
%and the definition of $\widehat{\mathfrak{F}}$ 
we obtain: 
%the following result.
%\color{blue} I'd mention that it is countable. 

\begin{lem}{\label{partial2}}
    $\widehat{\mathfrak{F}}$ is a countable partially ordered \textbf{MS4B}-frame.
\end{lem}

To obtain a model based on $\widehat{\mathfrak{F}}$, we define a valuation $\widehat{v}$ on $\widehat{\mathfrak{F}}$ by setting 
$$\widehat{v}(p)=\{\widehat{t}\in\widehat{X}\mid t\in v(p)\}$$ for $p\in S$ and $v(p)=\varnothing$ otherwise. As in \cref{subsec: MGrz}, %remarked earlier, 
we abbreviate satisfaction $\widehat{\mathfrak{M}},\widehat{t}\models\psi$ in the model $\widehat{\mathfrak{M}}=(\widehat{\mathfrak{F}},\widehat{v})$ by $\widehat{t}\models\psi$. %Consequently, we have the following result.

\begin{lem}[Truth Lemma]{\label{TL}}
    For any $\widehat{t}\in \widehat{X}$ and $\psi\in S$,
    \begin{equation*}
        \widehat{t}\models\psi \iff t\models\psi.
    \end{equation*}
\end{lem}
\begin{proof}
    The proof is by induction on the complexity of $\psi$. The base step %of the induction 
    is clear by the definition of $\widehat{v}$. 
    %Further, 
    By the induction hypothesis, the claim is clearly true when $\psi$ is a Boolean combination of subformulae. We only need to show the above when $\psi=\exists\delta$ or $\psi=\Diamond\delta$.

    Suppose $t\models\exists\delta$. If $t\models\delta$, then  $\widehat{t}\models\delta$ by the induction hypothesis, and so $\widehat{t}\mathrel{\widehat{E}}\widehat{t}$ implies $\widehat{t}\models\exists\delta$. Suppose $t\not\models\delta$. Then $\exists\delta\in W_t^\exists$, and there is a least $k<\omega$ for which  $\widehat{t}\in X_k$. Hence, by construction, there is 
    $\widehat{z}\in X_{k+1}$ such that $\widehat{t}\mathrel{E_{k+1}}\widehat{z}$ and $z\models\delta$. By the induction hypothesis, %we have 
$\widehat{z}\models\delta$. Since $E_{k+1}\subseteq\widehat{E}$, we have $\widehat{t}\mathrel{\widehat{E}}\widehat{z}$ and $\widehat{z}\models\delta$. Thus, $\widehat{t}\models\exists\delta$.

Conversely, suppose $\widehat{t}\models\exists\delta$. Then there is $\widehat{z}\in\widehat{X}$ such that $\widehat{t}\mathrel{\widehat{E}}\widehat{z}$ and $\widehat{z}\models\delta$. Take the least $k<\omega$ such that $\widehat{t}\mathrel{E_k}\widehat{z}$. By \crefdefpart{rmk}{rmk1}, $y\mathrel{E}z$, and by the induction hypothesis, $z\models\delta$. Hence, $y\models\exists\delta$.

    Next suppose $t\models\Diamond\delta$. %Once again, if we have 
    If $t\models\delta$, then  $\widehat{t}\models\delta$ by the induction hypothesis, implying that $\widehat{t}\models\Diamond\delta$ since $\widehat{t}\mathrel{\widehat{R}}\widehat{t}$. Suppose $t\not\models\delta$. Then $\Diamond\delta\in W_t^\Diamond$, and there is a least $k<\omega$ such that $\widehat{t}\in X_k$. By construction, there is $\widehat{z}\in X_{k+1}$ such that $\widehat{t}\mathrel{R_{k+1}}\widehat{z}$ and $z\models\delta$. The induction hypothesis implies that $\widehat{z}\models\delta$. Since $R_{k+1}\subseteq \widehat{R}$, we have $\widehat{t}\mathrel{\widehat{R}}\widehat{z}$ and $\widehat{z}\models\delta$. Thus, $\widehat{t}\models\Diamond\delta$.
    
    Conversely, suppose $\widehat{t}\models\Diamond\delta$. Then there is $\widehat{z}\in \widehat{X}$ such that $\widehat{t}\mathrel{\widehat{R}}\widehat{z}$ and $\widehat{z}\models\delta$. Choose the least $k, l<\omega$ such that $\widehat{t}\mathrel{R_k^l}\widehat{z}$. Because we take the reflexive-transitive closure at each step, we have a chain $$\widehat{t}=\widehat{x}_1\mathrel{R_k^l}\widehat{x}_2\mathrel{R_k^l}\dots\mathrel{R_k^l}\widehat{x}_n\mathrel{R_k^l}\widehat{x}_{n+1}=\widehat{z}.$$ 
    %where the relation $(\widehat{x}_j,\widehat{x}_{j+1})$ is introduced to $R_k$ during the course of the construction. 
    To see that $t\models\Diamond\delta$, we first show the following.
    \begin{clm}
        For any $j\geq2$, we have  $x_j\models\Diamond\delta\implies x_{j-1}\models \Diamond\delta$.
    \end{clm}
    \begin{proof}[Proof of the claim]
         %\color{blue} (should this be in the claim environment?) \color{black} 
         Suppose $x_j\models\Diamond\delta$. If $x_{j-1}\mathrel{R}x_j$, then %clearly 
    $x_{j-1}\models\Diamond\Diamond\delta$, which 
    %by transitivity 
    yields $x_{j-1}\models\Diamond\delta$. Suppose $x_{j-1}\mathrel{\cancel{R}}x_j$. By \crefdefpart{rmk}{rmk2}, $x_{j-1}\mathrel{Q}x_j$. Since $x_{j-1}\mathrel{Q}x_j$, $x_{j-1}\mathrel{\cancel{R}}x_j$, and $\widehat{x}_{j-1}\mathrel{R_k^l}\widehat{x}_j$, a proper $Q$-arrow was turned into a $R_k^l$-arrow. %during the construction. 
    %\color{blue} (this needs explanation; probably better to add a remark already in the previous subsection). \color{black} 
    This happens in only two stages of the construction: either in the vertical step to introduce $\Diamond$-witnesses or in the commutativity step (\hyperref[RC]{RC}) (see \Cref{Q-arrow1,Q-arrow2}). 
    %\color{blue} (again, add some detail; what about (LC)?). \color{black} 
    We handle each case separately.
    \begin{enumerate}
        \item Suppose $\widehat{x}_{j-1}\mathrel{R_k^l}\widehat{x}_j$ was introduced to add $\widehat{x}_j$ as a vertical $\Diamond$-witness for $\widehat{x}_{j-1}$ with respect to some $\Diamond\psi\in S$. Then by \labelcref{vw}, we must have $x_j\in\textbf{smax}\hspace{2pt}A$, where
        \begin{equation*}
            A = v(\Diamond\psi)\cap\bigcap\bigl\{ v(\neg\Diamond\alpha)\mid \Diamond\alpha\in S,\,  x_{j-1}\not\models \Diamond\alpha \bigr\}.
        \end{equation*}
        Hence, if $x_{j-1}\not\models\Diamond\delta$, we must have $x_j\not\models\Diamond\delta$, which contradicts our assumption that $x_j\models\Diamond\delta$. Therefore, $x_{j-1}\models\Diamond\delta$. %in this case.

        \item Suppose $\widehat{x}_{j-1}\mathrel{R_k^l}\widehat{x}_j$ was introduced to close under right commutativity. Then  (\hyperref[RC]{RC}) together with \Cref{r-comm} implies that 
        %we must have 
        $x_{j-1}\in V$, where 
        \begin{equation*}
            V=\bigcap\{v(\Diamond\psi)\mid \Diamond\psi\in S, x_j\models\Diamond\psi\}.
        \end{equation*}
        Thus, since $x_j\models\Diamond\delta$, we have $x_{j-1}\models\Diamond\delta$.\qedhere
    \end{enumerate}
    \end{proof}

    Now, applying \crefdefpart{rmk}{rmk4} to $\widehat{x}_{n+1}=\widehat{z}$ yields $x_{n+1}=z$, and by induction hypothesis we have $\widehat{z}\models\delta$. Thus, $x_{n+1}\models\delta$, and by the repeated application of $x_j\models\Diamond\delta\implies x_{j-1}\models\Diamond\delta$, we obtain $x_1\models\Diamond\delta$. From $\widehat{x}_1=\widehat{t}$ it follows that $x_1=t$ by \crefdefpart{rmk}{rmk4}, 
    %We obtain  from , 
    and hence $t \models\Diamond\delta$.
\end{proof}

We next show that $\widehat{\mathfrak{F}}=(\widehat{X},\widehat{R},\widehat{E})$ has bounded depth. 
For this we require working with the quotient of $\widehat{\mathfrak{F}}$
%a frame $(X,R,E)$ 
by the equivalence relation $\widehat{E}$. 
%we recall the notion of an \emph{$E$-skeleton}. 
This generalizes the well-studied notion of a skeleton 
%is well studied 
in modal logic (see, e.g., \cite[p.~68]{CZ}). 
%In particular, 
The skeleta of this form were first 
considered in the setting of monadic intuitionistic logics 
%by Suzuki 
(see, e.g., \cite{szk,BV}), and were  appropriately redefined for mm-logics in \cite{GBm,BM,GLf}.  

%\color{blue} Add some references. In the next definition, we shouldn't only reference our paper. \color{black}

\begin{defn} \label{E-skeleton}
    The {\em $E$-skeleton} of an \textbf{MK}-frame $\mathfrak{F}=(X,R,E)$ is the frame $\mathfrak{F}_0=(X_0,R_0)$, where $X_0$ is the quotient of $X$ by $E$. We denote elements of $X_0$ by $[x]$, where $x\in X$, and define
    \begin{equation*}
        [x]\mathrel{R_0}[y]\iff x\mathrel{Q}y.
    \end{equation*}
\end{defn}

Let $\widehat{\mathfrak{F}}_0=(\widehat{X}_0, \widehat{R}_0)$ be the $\widehat{E}$-skeleton of $\widehat{\mathfrak{F}}=(\widehat{X},\widehat{R},\widehat{E})$. The same argument as in \cite[Lem.~7.14]{BM} gives that $\widehat{\mathfrak{F}}_0$ is a poset. Indeed, if % which can be shown using the . We sketch this briefly. Suppose 
$[\widehat{x}]\mathrel{\widehat{R}_0}[\widehat{y}]\mathrel{\widehat{R}_0}[\widehat{x}]$, then 
%we can show that 
$x\mathrel{Q}y\mathrel{Q}x$, and by the strong maximality of $x$, we conclude that $x\mathrel{E}y$ (see \crefdefpart{smax_thms}{passive_Q}). Therefore, $\widehat{x}\mathrel{\widehat{E}}\widehat{y}$, which means that $[\widehat{x}]=[\widehat{y}]$. 
%\color{blue} (explain why; just referencing \cite[Lem.~7.11]{BM} isn't enough since the construction has changed). 

To see that $\widehat{\mathfrak{F}}$ has bounded depth, observe that the projection map $\pi\colon\widehat{\mathfrak{F}}\to\widehat{\mathfrak{F}}_0$, given by $\widehat{x}\mapsto [\widehat{x}]$, collapses each chain $C$ in $\widehat{\mathfrak{F}}$ to a chain $C_0$ in $\widehat{\mathfrak{F}}_0$. Thus, to see that the depth of each chain $C$ is bounded, it is enough to show that $\widehat{\mathfrak{F}}_0$ has bounded depth, 
%(see \Cref{chain1} below) 
and that the size of each set $\pi^{-1}[\widehat{y}]\cap C$ is uniformly bounded for each $[\widehat{y}]$ in $C_0$. 
%(see \Cref{chain2}) to conclude that  
%we can show that 
%the size of $C$ is bounded as well (see \Cref{bounded_depth}). 
%Therefore, $\widehat{\mathfrak{F}}$ has bounded depth.  
%\color{blue} This is an important passage, so needs better explanation. \color{black}
%Hence, we prove the following results. 

%\color{blue} The next lemma requires a proof since the construction has changed. Just referencing our previous paper isn't enough. Need to point out that the skeleton doesn't change. \color{black}

\begin{lem}{\label{chain1}}
    If $\widehat{C}$ is an $\widehat{R}_0$-chain \[[\widehat{y}_0]\mathrel{\widehat{R}_0}[\widehat{y}_1]\dots [\widehat{y}_{n-1}]\mathrel{\widehat{R}_0}[\widehat{y}_n]\dots\]  in $\widehat{\mathfrak{F}}_0=(\widehat{X}_0, \widehat{R}_0)$, then the length of $\widehat{C}$ is bounded by $2\cdot2^{|S|}$.
\end{lem}

\begin{proof}
    The proof 
    %of the above 
    follows that of \cite[Lem.~7.12]{BM} verbatim since the skeleton of $\widehat{\mathfrak{F}}$ is determined by the $\exists$ and $\Diamond$ steps of the construction, which do not require any modification.  
\end{proof}

\begin{lem}{\label{chain2}}
    If $\widehat{C}$ is an $\widehat{R}$-chain \[\widehat{x}_0\mathrel{\widehat{R}}\widehat{x}_1\dots \widehat{x}_{n-1}\mathrel{\widehat{R}}\widehat{x}_n\dots\] in $\widehat{X}$ such that $\widehat{x}_i\mathrel{\widehat{E}}\widehat{x}_j$ for all $i$, $j$, then the length of $\widehat{C}$ is bounded by $2^{|S|}$.
\end{lem}
\begin{proof}
    The proof is similar to \cite[Lem.~7.11]{BM}, but we have to account for an additional case due to the right commutativity step. First, observe that it is enough to establish the above claim when $\widehat{C}$ is a chain where every relation $\widehat{x}_i\mathrel{\widehat{R}}\widehat{x}_{i+1}$ is introduced at some stage of the construction (otherwise, we can decompose $\widehat{x}_i\mathrel{\widehat{R}}\widehat{x}_{i+1}$ further into a chain $\widehat{x}_i\mathrel{\widehat{R}}\widehat{y}_1\dots\widehat{y}_n\mathrel{\widehat{R}}\widehat{x}_{i+1}$, where the latter relations are indeed introduced through the course of the construction). 
    Now, $\widehat{x}_i\mathrel{\widehat{E}}\widehat{x}_j$ and $\widehat{x}_i\mathrel{\widehat{R}}\widehat{x}_j$ together imply that $\widehat{x}_i\mathrel{E_k}\widehat{x}_j$ and $\widehat{x}_i\mathrel{R_k}\widehat{x}_j$ for some $k<\omega$. By \crefdefpart{rmk}{rmk3}, %we get 
    $x_i\mathrel{E}x_j$ and $x_i\mathrel{R}x_j$. Therefore, we obtain the chain $C$  
    \begin{equation*}
        x_0\mathrel{R}x_1\dots x_{n-1}\mathrel{R}x_n\mathrel{R}\dots
    \end{equation*}
    in the canonical frame of \textbf{MGrzB}. 

    \begin{clm} \label{claim about sim}
     $x_i\mathrel{\cancel{\sim}_S} x_j$ whenever $i\neq j$. 
    \end{clm}

    \begin{proof}[Proof of the claim]
    %\color{blue} Should this claim be in the claim environment? \color{black} 
    Fix $i$ and without loss of generality assume that $i<j$. Define
    \begin{equation*}
         U= \bigcap\{v(\alpha)\mid \alpha\in S,\,  x_i\models\alpha\} \cap \bigcap \{v(\neg\beta) \mid \beta\in S,\,x_i\not\models\beta\}.
    \end{equation*}

    %Then we make the following 
    We observe that
    \begin{enumerate}[label={(\arabic*)}]
        \item $U$ is a definable subset of $X$ and $x_i\in U$. \label{obs1}

        \item For any $k$, we have $x_i\mathrel{\sim_S}x_k$ iff $x_k\in U$. \label{obs2}

        \item We have $\widehat{x}_i\mathrel{R_k}\widehat{x}_{i+1}$ and $\widehat{x}_i\mathrel{E_k}\widehat{x}_{i+1}$ for some $k<\omega$. By construction, such relations are introduced only in the horizontal $\Diamond$ step (see \Cref{horz_rmk}). Hence, by \labelcref{hw}, 
        %it follows that 
        there is 
        %a formula 
        $\Diamond\delta\in W_{x_i}^\Diamond$ such that $\widehat{x}_{i+1}$ is a $\Diamond$-witness for $\widehat{x}_i$ with respect to 
        %the formula 
        $\Diamond\delta$.  
        In particular, this implies that $x_i\models\Diamond\delta$, $x_i\not\models\delta$,  and $x_{i+1}\models\delta$. This 
        %clearly 
        shows that $x_{i+1}\notin U$, and hence that $x_i\mathrel{\cancel{\sim}}_Sx_{i+1}$. \label{obs3}
    \end{enumerate}

    We use the above to show that $x_i\mathrel{\cancel{\sim}_S} x_j$. 
    Since $\widehat{x}_i$ is introduced in either the $\exists$ step, $\Diamond$ step, or one of the commutativity steps, we address each separately.
    \begin{itemize}
        \item Suppose $\widehat{x}_i$ is introduced in the $\exists$ step. 
        %at some stage of the construction. 
        By \Cref{horz}, 
        %we must have 
$x_i\in\textbf{smax}\hspace{2pt}(E[V]\cap v(\psi))$ for some definable subset $V$ and $\exists\psi\in S$. Note that $U\subseteq v(\psi)$, and hence $E[V]\cap U\subseteq E[V]\cap v(\psi)$. This implies that $x_i\in\textbf{max}\hspace{2pt}(E[V]\cap U)$. We want to show that $x_i\mathrel{\cancel{\sim}_S} x_j$. Suppose otherwise. Then by putting $k=j$ in \labelcref{obs2} %\color{blue} (maybe label the above observations and add the label here to make it clear), \color{black} 
we get $x_j\in U$. Further, $x_j\mathrel{E}x_i$ and $x_i\in E[V]$ imply that $x_j\in E[V]$. Hence, $x_j\in E[V]\cap U$. However, %observe that 
$j>i+1$ since $x_{i+1}\notin U$.
Thus, $x_i\mathrel{R}x_{i+1}\mathrel{R}x_j$, where $x_i\in \textbf{max}\hspace{2pt}(E[V]\cap U)$, $x_{i+1}\notin E[V]\cap U$, and $x_j\in E[V]\cap U$. This contradicts that $\mathfrak{M}=(\mathfrak{F},v)$ is the canonical model of \textbf{MGrzB} (see the beginning of the section). 
%\color{blue} This should be mentioned earlier,  maybe in the theorem environment; the references could then go there and here it'd be enough to reference the theorem itself. \color{black}

        \item Suppose $\widehat{x}_i$ is introduced in the $\Diamond$ step. By \labelcref{hw} and \labelcref{vw}, we should have $x_i\in\textbf{max}\hspace{2pt} v(\psi)$ for some $\Diamond\psi\in S$. 
        %\color{blue} (what is $A$?). \color{black} 
        Because $U\subseteq v(\psi)$, and $x_i\in U$ by \labelcref{obs1}, we have $x_i\in\textbf{max}\hspace{2pt}U$. If $x_i\mathrel{\sim_S}x_j$, then $x_i\mathrel{R}x_{i+1}\mathrel{R}x_j$, where $x_i\in\textbf{max}\hspace{2pt}U$, $x_{i+1}\notin U$, and $x_j\in U$. Once again, this contradicts that $\mathfrak{M}$ is the canonical model of \textbf{MGrzB}.

        \item Suppose $\widehat{x}_i$ is introduced in the left commutativity step. The points used to close under left commutativity must be chosen using \Cref{l-comm}. 
    %\color{blue} (is this true?). \color{black} 
    Accordingly, we must have $x_i\in\textbf{smax}\hspace{2pt}E[V]$, where $V$ is a definable subset. Since $x_i \in U$ by \labelcref{obs1}, we have ${x}_i\in\textbf{max}\hspace{2pt}(E[V]\cap U)$, and the remainder of the proof proceeds in exactly the same way as in the first case.

        \item Finally, suppose $\widehat{x}_i$ is introduced in the right commutativity step. The points used to close under right commutatvity are chosen using \Cref{r-comm}. Therefore, we must have $x_i\in\textbf{smax}\hspace{2pt}(E[W]\cap V)$, where both $W$ and $V$ are some definable sets. Since $x_i \in U$ by \labelcref{obs1},
        %From \labelcref{obs1}, 
        we obtain $x_i\in\textbf{max}\hspace{2pt}(E[W]\cap V\cap U)$, and the rest of the argument again follows the first case.\qedhere
    \end{itemize}
    \end{proof}
    
    Therefore, the chain $C$, and thus $\widehat{C}$, can have at most $2^{|S|}$ points.
\end{proof}

Combining \Cref{chain1,chain2}, 
%\color{blue} (why is ``Lemmas'' added manually?) \color{black}, 
we obtain: %the following.

\begin{thm}\label{bounded_depth}
    The depth 
    %\color{blue} (are we supposed to say $\widehat{R}$-depth?) \color{black} 
    of $\widehat{\mathfrak{F}}$ is bounded by $2^{2|S|+1}$.
\end{thm}
\begin{proof}
    Let $$ \widehat{x}_0\mathrel{\widehat{R}}\widehat{x}_1\dots\widehat{x}_{n-1}\mathrel{\widehat{R}}\widehat{x}_n\mathrel{\widehat{R}}\dots$$ be a chain $C$ in $\widehat{\mathfrak{F}}$, and let $C$ be indexed by $I\subseteq\omega$. Define $f\colon I\to I$ by $$f(i)=\min\{k\in I\mid \widehat{x}_i\mathrel{\widehat{E}}\widehat{x}_k\}.$$ Let $s_0,s_1,\dots, s_m,\dots$ be an enumeration of $f(I)$ such that $s_i<s_j$ whenever $i<j$. 
    %Then we obtain a chain
    Consider the corresponding chain $$[\widehat{x}_{s_0}]\mathrel{\widehat{R}_0}[\widehat{x}_{s_1}]\dots [\widehat{x}_{s_{m-1}}]\mathrel{\widehat{R}_0}[\widehat{x}_{s_m}]\dots$$ in $\widehat{\mathfrak{F}}_0$. By \Cref{chain1}, %$C_0$ 
    this chain (and hence $f(I)$) has size at most $2\cdot 2^{|S|}$. For a fixed $m$, consider the chain $C\cap [\widehat{x}_{s_m}]$. By \Cref{chain2}, 
    %we have 
    $|C\cap [\widehat{x}_{s_m}]|\leq 2^{|S|}$. Now observe that
    \begin{equation*}
        C=\bigsqcup_{m\in f(I)}(C\cap [\widehat{x}_{s_m}]).
    \end{equation*}
    Thus, $|C|\leq |f(I)|\cdot 2^{|S|}\leq 2^{2|S|+1}$, as required.
\end{proof}

%\color{blue} Maybe add a short proof explaining how the upper bound is obtained. \color{black}

\section{FMP and axiomatization}\label{sec5}

%\color{blue}
%Add a paragraph briefly explaining the content of the section.
%\color{black}

In this final section, we use the selective filtration technique developed in the previous section to show that 
%consolidate all of our observations and show that the systems 
the mm-logics $\textbf{MGrzB}$, $\textbf{M}^+\textbf{GrzB}$, and $\textbf{MGLB}$ have the fmp. We then use the fmp of the latter two 
%logics 
to show that $\textbf{GL}\times\textbf{S5}$ and $\textbf{Grz}\times\textbf{S5}$ are finitely axiomatizable. This resolves the question of finite axiomatization for two of the logics from the list given in \cite[p.~137]{gab98}. %\color{blue} Mention that this resolves the open problem mentioned in Gabbay-Shehtman. \color{black} 

%\subsection{FMP for \textbf{MGrzB}, \texorpdfstring{$\textbf{M}^+\textbf{GrzB}$}{M+GrzB}, and \textbf{MGLB}}

    %As mentioned above, we will prove that these systems have the fmp. 
    The key feature of our construction is that, for each ${\varphi\notin\textbf{MGrzB}}$, it produces a countable \textbf{MGrzB}-frame of bounded depth refuting $\varphi$. We then use %the local tabularity of $\textbf{MGrz}[n]$ 
    %and $\textbf{MGL}[n]$ 
    %(see 
    \Cref{locally_finite} to produce a finite \textbf{MGrzB}-frame refuting $\varphi$, thus yielding the fmp. 

%\color{blue}
%Start by saying what you are after, then have usual theorem/proof environment. Same applies to the next theorem.
%\color{black}

\begin{thm}\label{MGrzB_fmp}
    \textbf{MGrzB} has the fmp.
\end{thm}

\begin{proof}
    Let $\varphi$ be an $\mathcal{L}_{\Diamond\exists}$-formula such that $\varphi\notin\textbf{MGrzB}$. Then $\varphi$ 
    %formula 
    is refuted in the canonical model $\mathfrak{M}=(\mathfrak{F},v)$ for $\textbf{MGrzB}$. Applying our construction to $\mathfrak{M}$ produces a countable model  $\widehat{\mathfrak{M}}=(\widehat{\mathfrak{F}},\widehat{v})$ such that 
    %can be constructed. From \Cref{TL}, we observe that 
    $\widehat{\mathfrak{F}}\not\models\varphi$ %while \Cref{bounded_depth} shows that 
    and $\widehat{\mathfrak{F}}$ is an  $\textbf{MGrzB}[n]$-frame, where $n=2^{2|S|+1}$ and $S$ is the set of subformulae of $\varphi$. 
    %Therefore, $\varphi\notin\textbf{MGrz}[n]\textbf{B}$.  
    Since $\textbf{MGrzB}[n]$ is locally tabular (see \Cref{locally_finite}), there is a finite p-morphic image 
    %finite $\textbf{MGrzB}[n]$-frame 
    $\mathfrak{G}$ of $\widehat{\mathfrak{F}}$ such that $\mathfrak{G}\not\models\varphi$ (see \Cref{p-morphicimage}).  
    %\color{blue} (this can be explained in \cref{sec3}). \color{black} 
    Since $\mathfrak{G}$ is an \textbf{MGrzB}-frame, and the result follows.
\end{proof}

%\color{red}
%Done
%\color{blue}
%I removed the proof. But another option is to have the above paragraph in the proof environment after the theorem.  
%\color{black}

\begin{thm}\label{fmp}
Both $\textbf{M}^+\textbf{GrzB}$ and $\textbf{MGLB}$ have the fmp.
\end{thm}
\begin{proof}
    The fmp for both $\textbf{M}^+\textbf{GrzB}$ and \textbf{MGLB} can be established by making modifications to the construction for \textbf{MGrzB}. For $\textbf{M}^+\textbf{GrzB}$, these modifications simply involve the omission of the horizontal $\Diamond$ step since the $E$-clusters of maximal points in the canonical model for $\textbf{M}^+\textbf{GrzB}$ are clean (recall \Cref{clean} and \cite[Lem.~4.8]{GBm}). For $\textbf{MGLB}$, we work with the strong irreflexive maxima of definable subsets. The $E$-clusters of such points are clean, which again allows us to omit the horizontal $\Diamond$ step. For both of these systems, a version of \Cref{TL} applies to the construction, yielding that the depth of the constructed frame $\widehat{\mathfrak{F}}$ is at most $2^{|S|+1}$. Note that the upper bound for the depth of $\widehat{\mathfrak{F}}$ is lower here since the estimation of the depth of the frame $\widehat{\mathfrak{F}}$ reduces to knowing the depth of its $E$-skeleton $\widehat{\mathfrak{F}}_0$, which can be deduced using \Cref{chain1}. 
    %\color{red} done \color{blue} This also applies to  $\textbf{M}^+\textbf{GrzB}$. \color{black}
\end{proof}

\begin{rmk}
%\color{red}
It is known that $\textbf{MGL}$ and $\textbf{M}^+\textbf{Grz}$ have the fmp (see \cite{japaridze88, japaridze, GBm, BM}).
%were shown to possess the fmp \color{blue} (the fmp of $\textbf{MGL}$ goes back to Japaridze), \color{black} 
%and these results were established using selective filtration. 
The selection procedure in the absence of the monadic Barcan formula is simpler since it is not necessary to revisit $E$-clusters that have been populated with enough witnesses. Hence, the size of $E$-clusters is finite. In contrast, with the construction proposed here, the requirement of right commutativity forces updates on previously constructed $E$-clusters, which may result in  infinite $E$-clusters. 

We emphasize that the constructions proposed in \cite{japaridze88,japaridze,GBm} produce an expanding relativized product frame. This is not the case with the selection technique we develop here; that is, the countermodel produced by the construction is \emph{not} a product frame. Indeed, our construction does not prevent the existence of $x$ and $y \ne z$ such that $x\mathrel{R}y$, $x\mathrel{R}z$, and $y\mathrel{E}z$. Thus, the resulting frame is not necessarily a product frame. %$\mathfrak{F}\times\mathfrak{G}=(X\times Y,R,E)$, given $a,b,c\in X\times Y$, if $a\mathrel{R}b$, $a\mathrel{R}c$, and $b\mathrel{E}c$, then $b=c$. However, this property need not hold in the frame constructed using our method. \color{blue} Briefly explain why. \color{black}
%Indeed, the question of whether these logics have the product fmp (see, e.g., \cite[p.~132]{mdim}) remains open. 
%\color{blue}
%There's another thing to point out: the construction given say in \cite{GBm} does produce an expanding relativized product frame, while the construction here doesn't. Need to explain. 
%Discuss how the above compares to how fmp was proved for these logics without Barcan.
%\color{black}
\end{rmk}

%\color{red}
%\begin{proof}
%  The proof is analogous to that of \Cref{MGrzB_fmp}. In particular, our construction can be used to show that for any non-theorem $\varphi$, there is an appropriate frame of finite depth that refutes $\varphi$. Hence, by \Cref{locally_finite}, both the above claims follow.
%\end{proof}
%\color{black}

Since the above three logics are %$\textbf{MGrzB}$, $\textbf{M}^+\textbf{GrzB}$, and $\textbf{MGLB}$ are 
finitely axiomatizable, as an immediate consequence we obtain: 
%logics with the fmp, we obtain the following as a consequence of Harrop's theorem (see, e.g., \cite[Thm.~16.3]{CZ}).

\begin{thm}\label{decidable}
    The logics $\textbf{MGrzB}$, $\textbf{M}^+\textbf{GrzB}$, and $\textbf{MGLB}$ are decidable.
    %\begin{enumerate}[label={(\arabic*)}]
    %    \item[] 
    %    \item $\textbf{MGrzB}$ is decidable. \label{dec1}
    %    \item $\textbf{M}^+\textbf{GrzB}$ is decidable. \label{dec2}
    %    \item $\textbf{MGLB}$ is decidable. \label{dec3}
    %\end{enumerate}
\end{thm}

%\subsection{Axiomatizing \texorpdfstring{$\textbf{Grz} \times \textbf{S5}$}{Grz x S5} and \texorpdfstring{$\textbf{GL} \times \textbf{S5}$}{GL x S5}}

We now use the fmp of the latter two logics %above results 
to obtain finite axiomatizations for $\textbf{GL} \times \textbf{S5}$ and $\textbf{Grz} \times \textbf{S5}$. In particular, we show that $\textbf{GL}\times\textbf{S5}$ is product matching, and that $\textbf{Grz}\times\textbf{S5}=\textbf{M}^+\textbf{GrzB}$.

%\color{red} done
%\color{blue}
%It might make sense to combine  \Cref{product_matching_GL,axiomatization} into one theorem with two items. 
%\color{black}

\begin{thm}\label{main}
    \begin{enumerate}[label={(\arabic*)}]
        \item[] 
        \item  $\textbf{GL}\times\textbf{S5}$ is product matching. \label{product_matching_GL}

        \item $\textbf{Grz}\times\textbf{S5}=\textbf{M}^+\textbf{GrzB}$. \label{axiomatization}
    \end{enumerate}
\end{thm}
\begin{proof}
    \labelcref{product_matching_GL} The inclusion $\textbf{MGLB}\subseteq \textbf{GL}\times\textbf{S5}$ always holds (see, e.g., \cite[Prop.~3.8]{mdim}). For the other inclusion, suppose $\varphi\notin\textbf{MGLB}$. By \Cref{fmp}, there is a finite $\textbf{MGLB}$-frame $\mathfrak{F}$ such that $\mathfrak{F}\not\models\varphi$. Let $n=d(\mathfrak{F})$. Then $\mathfrak{F}$ is an $\textbf{MGLB}[n]$-frame, and thus $\varphi\notin\textbf{MGLB}[n]$. By \Cref{matching_GLn}, $\varphi\notin \textbf{GL}[n]\times\textbf{S5}$. Since $\textbf{GL}\times\textbf{S5}\subseteq\textbf{GL}[n]\times\textbf{S5}$, we conclude that $\varphi\notin\textbf{GL}\times\textbf{S5}$.
    
    \labelcref{axiomatization} By \Cref{casari_rmk2},  $\textbf{M}^+\textbf{GrzB}\subseteq\textbf{Grz}\times\textbf{S5}$. For the other inclusion, %Hence, to show $\textbf{Grz}\times\textbf{S5}\subseteq\textbf{M}^+\textbf{GrzB}$, 
    suppose $\varphi\notin\textbf{M}^+\textbf{GrzB}$. By \Cref{fmp}, there is a finite $\textbf{M}^+\textbf{GrzB}$-frame $\mathfrak{F}$ such that $\mathfrak{F}\not\models\varphi$. Let $n=d(\mathfrak{F})$. Then $\mathfrak{F}$ is an $\textbf{M}^+\textbf{GrzB}[n]$-frame, and hence $\varphi\notin\textbf{M}^+\textbf{GrzB}[n]$. 
    %However, we know that $\textbf{M}^+\textbf{GrzB}[n]=\textbf{Grz}[n]\times\textbf{S5}$ from 
    By \Cref{axiomatization_n}
    %. Thus, we have 
    $\varphi\notin\textbf{Grz}[n]\times\textbf{S5}$. Thus,  $\varphi\notin\textbf{Grz}\times\textbf{S5}$ since $\textbf{Grz}\times\textbf{S5}\subseteq\textbf{Grz}[n]\times\textbf{S5}$. 
\end{proof}

%\color{red} Maybe we need a different subsection for this. \color{black}

%\color{blue}
%Use fmp and stuff from section 3. Remove 5.6 and add short proof to 5.7.
%\color{black}

%\color{red} Removed 5.6. Proof added to the translation theorem. \color{black}

Recall that the translation of $\textbf{Grz}$ into $\textbf{GL}$ lifts to yield $\varphi\in\textbf{M}^+\textbf{Grz}\iff\varphi^+\in\textbf{MGL}$ \cite[Thm.~4.12]{GBm}. We next show that the addition of the monadic Barcan formula preserves the faithfulness of this translation. 
%In light of the completeness of $\textbf{MGLB}$ and $\textbf{M}^+\textbf{GrzB}$, we obtain the following result.

\begin{thm}\label{embedding}
     For each $\mathcal{L}_{\Diamond\exists}$-formula $\varphi$,  $$\varphi\in\textbf{M}^+\textbf{GrzB}\iff \varphi^+\in\textbf{MGLB}.$$
\end{thm}

\begin{proof}
    Suppose $\varphi\notin\textbf{M}^+\textbf{GrzB}$. By \Cref{fmp}, there is a finite $\textbf{M}^+\textbf{GrzB}$-frame $\mathfrak{F}$ such that $\mathfrak{F}\not\models\varphi$. Since $\mathfrak{F}$ is also an $\textbf{M}^+\textbf{GrzB}[n]$-frame for some $n<\omega$, we have $\varphi\notin\textbf{M}^+\textbf{GrzB}[n]$. By \crefdefpart{translation}{translation4}, $\varphi^+\notin\textbf{MGLB}[n]$. Since $\textbf{MGLB}\subseteq\textbf{MGLB}[n]$, we conclude that $\varphi^+\notin\textbf{MGLB}$.

    Conversely, let $\varphi^+\notin\textbf{MGLB}$. By \Cref{fmp}, there is a finite $\textbf{MGLB}$-frame $\mathfrak{F}$ such that $\mathfrak{F}\not\models\varphi^+$. Because $\mathfrak{F}$ is a frame for $\textbf{MGLB}[n]$ for some $n<\omega$, we deduce that $\varphi^+\notin\textbf{MGLB}[n]$. By \crefdefpart{translation}{translation2}, $\varphi\notin\textbf{M}^+\textbf{GrzB}[n]$. Therefore, $\varphi\notin\textbf{M}^+\textbf{GrzB}$ since $\textbf{M}^+\textbf{GrzB}\subseteq\textbf{M}^+\textbf{GrzB}[n]$.
\end{proof}

Finally, we use our criterion from \cite[Thm.~5.4]{BM} to show that $\textbf{M}^+\textbf{GrzB}$ and $\textbf{MGLB}$ axiomatize the one-variable fragments of $\textbf{Q}^+\textbf{GrzB}$ and $\textbf{QGLB}$, respectively. Recall that there is a Wajsberg-style translation that associates with each $\mathcal{L}_{\Diamond\exists}$-formula $\varphi$ a monadic predicate formula $\varphi^t$ (see, e.g., \cite[Def.~13.4]{gab98}, \cite[p.~143]{mdim}, or \cite[p.~9]{BM}). 
For an mm-logic $\textbf{M}$ and a predicate modal logic $\textbf{Q}$, we say that $\textbf{M}$ {\em axiomatizes the one-variable fragment of} $\textbf{Q}$ if for any $\mathcal{L}_{\Diamond\exists}$-formula $\varphi$, we have $$\varphi\in \textbf{M}\iff\varphi^t\in \textbf{Q}.$$
This translation has a semantic counterpart, for which we recall that a {\em Kripke bundle} is a triple $\mathfrak{B}=(\mathfrak{F},\pi,\mathfrak{F}_0)$, where $\mathfrak{F}$ and $\mathfrak{F}_0$ are propositional Kripke frames, and $\pi\colon\mathfrak{F}\to\mathfrak{F}_0$ is an onto p-morphism (see, e.g., \cite[Def.~5.2.5]{GSS}). To see how predicate formulae are interpreted in a Kripke bundle, consult \cite[Def.~5.2.8]{GSS}. There is a functor $\mathscr{B}$ \cite[Prop.~4.3]{BM} that associates with each $\textbf{MK}$-frame $\mathfrak{F}=(X,R,E)$ the Kripke bundle $\mathscr{B}(\mathfrak{F})=((X,R),\pi,(X_0,R_0))$, where $(X_0,R_0)$ is the $E$-skeleton of $\mathfrak{F}=(X,R,E)$ (see \Cref{E-skeleton}), and $\pi\colon (X,R)\to(X_0,R_0)$ is the quotient map $x\mapsto[x]$. In \cite[Thm.~5.3]{BM}, we showed that
\begin{equation}
    \mathfrak{F}\models\varphi \iff 
    \mathscr{B}(\mathfrak{F})\models \varphi^t. \label{bundle_translation}\tag{$*$}
\end{equation}
For a product frame $\mathfrak{F}\times\mathfrak{G}$, where %$\mathfrak{F}=(X,R)$ and 
$\mathfrak{G}=(D,D^2)$ is an $\textbf{S5}$-frame, $\mathscr{B}(\mathfrak{F}\times\mathfrak{G})$ can be identified with the predicate Kripke frame with constant domain $(\mathfrak{F},D)$ (see, e.g., \cite[p.~149]{mdim}), 
%\cite[Lem.~13.5]{gab98}), 
and
%\color{red}
%In particular, we point out that $\mathscr{B}$ associates with each rooted product frame $\mathfrak{F}\times\mathfrak{G}$ a Kripke bundle corresponding to a predicate Kripke frame \cite[Def.~3.2.2]{GSS} over $\mathfrak{F}$. Specifically, let $\mathfrak{F}=(X,R_X)$, $\mathfrak{G}=(D,D^2)$, and $\mathfrak{F}\times\mathfrak{G}=(X\times D,R,E)$; then the $E$-skeleton of $\mathfrak{F}\times\mathfrak{G}$ is isomorphic to $\mathfrak{F}$. Therefore, we may regard $\mathscr{B}(\mathfrak{F}\times\mathfrak{G})=((X\times D,R),\pi,(X,R_X))$, where the map $\pi$ is given by $(x,a)\mapsto x$. This in turn corresponds to the predicate Kripke frame $((X,R_X),D)$, and it is easy to verify that for any predicate formula $\varphi$, $\mathscr{B}(\mathfrak{F}\times\mathfrak{G})\models\varphi\iff ((X,R_X),D)\models\varphi$. Hence, for all intents and purposes, we can assume that $\mathscr{B}(\mathfrak{F}\times\mathfrak{G})$ is a predicate Kripke frame over $\mathfrak{F}$.
%\color{blue} Ok, but it should be pointed out that Ono-Suzuki work with posets, however the result holds in full generality. Is this not addressed in either Shehtman-Skvortsov or GSS? Also, have predicate Kripke frames been mentioned yet?
%Ono-Suzuki basically do this in their paper, but for intuitionistic Kripke frames. You can't just wave your hands, proofs need to be formal. 
%\color{red} 
our criterion from \cite[Thm.~5.4]{BM} 
%for the faithfulness of the translation can be 
reduces to the following: 

\begin{thm}\label{criterion}
Let $\textbf{L}$ be a Kripke complete propositional modal logic and $\textbf{Q}$ a predicate modal logic. Then $\textbf{L}\times\textbf{S5}$ axiomatizes the one-variable fragment of $\textbf{Q}$ provided
\begin{enumerate}[label={(\arabic*)}]
    \item $\varphi\in\textbf{L}\times\textbf{S5}\implies\varphi^t\in\textbf{Q}$ for any $\mathcal{L}_{\Diamond\exists}$-formula; \label{criterion1}
    \item $\textbf{Q}$ is sound with respect to the class $\{\mathscr{B}(\mathfrak{F}\times\mathfrak{G})\mid \mathfrak{F} \text{ is an $\textbf{L}$-frame and } \mathfrak{G}=(D,D^2) 
    %\text{ is a rooted }\textbf{L}\times\textbf{S5}\text{-frame}
    \}$. 
    %\color{blue} This is still not how you use it below! \color{black} 
    \label{criterion3}
\end{enumerate}
\end{thm}

\begin{proof}
    %As usual, we identify a logic with the set of its theorems. Thus, we show that $\varphi\in \textbf{L}\times\textbf{S5}\iff\varphi^t\in \textbf{Q}$. 
    Due to \labelcref{criterion1}, we only need to show that $\varphi\notin\textbf{L}\times\textbf{S5}\implies\varphi^t\notin\textbf{\textbf{Q}}$. Since $\varphi\notin\textbf{L}\times\textbf{S5}$, $\varphi$ is refuted on a product frame $\mathfrak{F}\times\mathfrak{G}$, where $\mathfrak{F}$ is an $\textbf{L}$-frame and $\mathfrak{G}=(D,D^2)$.
    %$\times\textbf{S5}$-frame $\mathfrak{F}\times\mathfrak{G}$. 
    By \labelcref{criterion3}, $\mathscr{B}(\mathfrak{F}\times\mathfrak{G})$ is a predicate Kripke frame for $\textbf{Q}$, and by \labelcref{bundle_translation}, $\mathscr{B}(\mathfrak{F}\times\mathfrak{G})\not\models\varphi^t$. Thus, $\varphi^t\notin\textbf{Q}$.
\end{proof}

\begin{rmk}
    Since we do not assume that $\textbf{Q}$ is Kripke complete, the above theorem generalizes \cite[Thm.~3.21]{mdim} (as well as \cite[Thm.~13.8(1)]{gab98} because the only restriction we place on $\mathbf{L}$ is that it is Kripke complete). 
    %In fact, those references only apply to those predicate logics in \Cref{fragments} below that are Kripke complete (note that $\mathbf{QGLB}$ is Kripke incomplete \color{blue} what do we know about Kripke incompleteness of $\textbf{Q}^+\textbf{GrzB}$?). \color{black} 
    %We also point out that existing results on axiomatizing one-variable fragments may assume completeness of $\textbf{Q}$ (e.g., \cite[Thm.~3.21]{mdim}), or place some restrictions on the logic $\textbf{L}$ (e.g., \cite[Thm.~13.8(1)]{gab98}). We don not place such restrictions, and as such our criterion generalizes the existing ones.
\end{rmk}

In light of \Cref{criterion},
%this criterion, 
we have the following: 
%\color{blue} what about $\textbf{MGrzB}$? \color{black}
%results.

\begin{thm}\label{fragments}
    \begin{enumerate}[label={(\arabic*)}]
        \item[] 
        \item $\textbf{M}^+\textbf{GrzB}$ axiomatizes the one-variable fragment of $\textbf{Q}^+\textbf{GrzB}$. \label{Q+GrzB}

        \item $\textbf{MGLB}$ axiomatizes the one-variable fragment of $\textbf{QGLB}$. \label{QGLB}

        \item $\textbf{M}^+\textbf{GrzB}[n]$ axiomatizes the one-variable fragment of $\textbf{Q}^+\textbf{GrzB}[n]$. \label{Q+GrzBn}

        \item $\textbf{MGLB}[n]$ axiomatizes the one-variable fragment of $\textbf{QGLB}[n]$. \label{QGLBn}
    \end{enumerate}
\end{thm}

\begin{proof}
We prove \labelcref{Q+GrzB}. The proofs for \labelcref{QGLB,Q+GrzBn,QGLBn} are similar. Since 
%we have 
$\textbf{Grz}\times\textbf{S5}=\textbf{M}^+\textbf{GrzB}$, we view $\textbf{M}^+\textbf{GrzB}$ as the product logic $\textbf{Grz}\times\textbf{S5}$. The first item of \Cref{criterion} is easily checked since the translations of all the axioms of $\textbf{M}^+\textbf{GrzB}$ are theorems of $\textbf{Q}^+\textbf{GrzB}$. For the second item, if $\mathfrak{F}\times\mathfrak{G}$ is a product frame with $\mathfrak{F}$ a $\textbf{Grz}$-frame and $\mathfrak{G}=(D,D^2)$,
%$\times\textbf{S5}$-frame, 
then $\mathscr{B}(\mathfrak{F}\times\mathfrak{G})$ is a predicate $\mathbf{Grz}$-frame with constant domain. 
%over $\mathfrak{F}$, which is a $\textbf{Grz}$-frame. 
Hence, $\mathscr{B}(\mathfrak{F}\times\mathfrak{G})$ is a $\textbf{QGrzB}$-frame. Since $\mathscr{B}(\mathfrak{F}\times\mathfrak{G})$ is a Noetherian predicate Kripke frame, the Casari formula is intuitionistically valid in it \cite[Thm.~3(2)]{ono}. The G\"{o}del translation of the Casari formula is therefore valid in $\mathscr{B}(\mathfrak{F}\times\mathfrak{G})$ (see, e.g., \cite[Prop.~3.2.25]{GSS}), and hence it is a $\textbf{Q}^+\textbf{GrzB}$-frame. Thus, \Cref{criterion} applies to yield the result. 
\end{proof}

\begin{rmk}
    In the above proof, the finite axiomatization of $\textbf{Grz}\times\textbf{S5}$ was used in a crucial way to verify the first item of \Cref{criterion}.
    %in the proof above. 
    In general, proving the implication $\varphi\in\textbf{L}\times\textbf{S5}\implies\varphi^t\in \textbf{Q}$ may not be a trivial matter.
    %hard if $\textbf{Q}$ is a non-trivial logic. 
\end{rmk}

\begin{rmk}
    %We also note that 
    %On the other hand, 
    In contrast, $\textbf{MGrzB}[n]$ is \emph{not} the one-variable fragment of any predicate modal logic. Indeed, the translations of the axioms of $\textbf{MGrzB}[n]$ result in a predicate modal logic containing $\textbf{QGrzB}[n]$. However, by \crefdefpart{fragments}{Q+GrzBn}, the one-variable fragment of $\textbf{QGrzB}[n]$ is axiomatized by $\textbf{M}^+\textbf{GrzB}[n]$ since $\textbf{QGrzB}[n]=\textbf{Q}^+\textbf{GrzB}[n]$. The latter equality follows from the fact that $\textbf{QGrzB}[n]$ is Kripke complete: the predicate extension of any universal propositional modal logic 
    %(see, e.g., \cite[Def.~1.12.5]{GSS}) 
    %in the presence of 
    together with the Barcan formula is Kripke complete (see, e.g., \cite[Thm.~7.4.7]{GSS}), and $\textbf{Grz}[n]$ is indeed a universal modal logic. Due to the Kripke completeness of $\textbf{QGrzB}[n]$, the G\"odel translation of the Casari formula is provable in it since, as we pointed out in the proof of \Cref{fragments}, 
    %the latter 
    it is valid in all Noetherian predicate Kripke frames. 
    %(see, e.g., \cite[Thm.~3(2)]{ono}). 

    It remains open whether $\textbf{MGrzB}$ axiomatizes the one-variable fragment of $\textbf{QGrzB}$. The difficulty of determining this lies 
    %is due to 
    in the fact that for an $\textbf{MGrzB}$-frame $\mathfrak{F}$, the set of formulae valid in the Kripke bundle $\mathscr{B}(\mathfrak{F})$ may not be closed under substitution (see, e.g., \cite[p.~350]{GSS}). In particular, substitution instances of the Barcan formula need not be valid in $\mathscr{B}(\mathfrak{F})$ even if the Barcan formula is (a simple counterexample can be found in \cite[p.~92]{suzuki93}). This warrants investigating whether we can find a class $\mathcal K$ of $\textbf{MGrzB}$-frames such that $\textbf{MGrzB}$ is complete with respect to $\mathcal K$ and for each $\mathfrak{F} \in \mathcal K$, all substitution instances of the Barcan formula are valid in $\mathscr{B}(\mathfrak{F})$. 
    %In this case, our criterion from \cite[Thm.~5.4]{BM} applies; however, it is unclear if such a class of $\textbf{MGrzB}$-frames exists. 
    %\color{blue} Bundles are used again, so it is important to give the reader at least some info about them. The business about validity of Barcan also requires references, otherwise the reader is left hanging. \color{black} 
    %\color{blue} Briefly explain the diffculty. Also, should we say ``is" or rather ``axiomatizes" the one-variable fragment? \color{red} I think it is better to say ``axiomatizes''. \color{black}
\end{rmk}

%\color{blue}
%Where will the business about monadic fragments of predicate logics be mentioned?
%\color{black}

\subsection*{Acknowledgements} 

We are thankful to Ilya Shapirovsky for useful conversations.  

\printbibliography
\end{document}